\documentclass{imamatgeneric}
\usepackage{amsmath,amssymb,color,graphicx,multirow}
\usepackage{placeins} % For \FloatBarrier command
\begin{document}

\title{Solution of Wiener-Hopf and Fredholm integral equations by fast Hilbert and Fourier transforms}
\shorttitle{Solution of Wiener-Hopf and Fredholm integral equations} %%%for recto running head
\shortauthorlist{G.~Germano, C.~E.~Phelan, D.~Marazzina, G.~Fusai}

\author{
\name{Guido Germano}
\address{Department of Computer Science, University College London \\
Systemic Risk Centre, London School of Economics and Political Science}
\email{g.germano@ucl.ac.uk}
\name{Carolyn E.\ Phelan}
\address{Department of Computer Science, University College London}
\email{carolyn.phelan.14@ucl.ac.uk (corresponding author)}
\name{Daniele Marazzina}
\address{Dipartimento di Matematica, Politecnico di Milano}
\email{daniele.marazzina@polimi.it}
\name{Gianluca Fusai}
\address{Business School (formerly Cass), City, University of London \\
Dipartimento di Studi per l'Economia e l'Impresa,
Universit\`a del Piemonte Orientale Amedeo Avogadro, Novara}
\email{gianluca.fusai.1@city.ac.uk, gianluca.fusai@unipmn.it}}

\maketitle

\begin{abstract}
{We present numerical methods based on the fast Fourier transform (FFT) to solve convolution integral equations on a semi-infinite interval (Wiener-Hopf equation) or on a finite interval (Fredholm equation). We extend and improve a FFT-based method for the Wiener-Hopf equation due to Henery, expressing it in terms of the Hilbert transform, and computing the latter in a more sophisticated way with sinc functions. We then generalise our method to the Fredholm equation reformulating it as two coupled Wiener-Hopf equations and solving them iteratively. We provide numerical tests and open-source code.}
{Wiener-Hopf, Fredholm, integral equation, fast Fourier transform, fast Hilbert transform.}
\end{abstract}

%\keywords{Wiener-Hopf, Fredholm, integral equation, fast Fourier transform, fast Hilbert transform.}

\section{Introduction}

We consider the linear integral equation of convolution type with constant limits of integration
\begin{equation}
\label{eq:WHF}
\lambda f(x) - \int_a^b k(x-x') f(x') dx' = g(x),\quad x \in (a,b),
\end{equation}
where $f(x)$ is the unknown function, $k(x)$ is a given kernel, and $g(x)$ is a given so-called forcing function. The domain of $f(x)$ and $g(x)$ is $(a,b)$, the domain of $k(x)$ is $(a-b,b-a)$; the endpoints can be included if they are finite. If $a = -\infty$ or $b = +\infty$ Eq.~(\ref{eq:WHF}) is called a Wiener-Hopf equation \citep{Wiener1931,Noble1958,Krein1963,Polyanin1998,Lawrie2007}; if both integration limits are finite, it is called a Fredholm equation \citep{Fredholm1903,Whittaker1927,Polyanin1998}. The latter case is also called a Wiener-Hopf equation on a finite interval \citep{Voronin2004} or, because of an application in electrotechnics, a longitudinally modified Wiener-Hopf equation (LMWHE), while the former case is also called a classical Wiener-Hopf equation (CWHE) \citep{Daniele2014}. If $\lambda = 0$ it is an equation of the first kind; if $\lambda \neq 0$ it is an equation of the second kind. In the latter case it can be assumed that $\lambda = 1$, dividing the kernel and the forcing function by values of this parameter different from 1. Historically these equations arose in physics, e.g.\ to describe diffraction in the presence of an impenetrable wedge or of planar waveguides \citep{Daniele2007}, but also for problems in crystal growth, fracture mechanics, flow mechanics \citep{Choi2005}, geophysics, and diffusion \citep{Lawrie2007}. The connection of the Wiener-Hopf equation with probabilistic problems was noticed by \cite{Spitzer1957} and is discussed by \cite{Feller1971}, together with the application of Fourier transform methods to stochastic processes. More recently these equations have become of interest in finance to price discretely monitored path-dependent options like barrier, first-touch, lookback (or hindsight), quantile and Bermudan options \citep{Fusai2006,Green2010,Fusai2012,Marazzina2012,Fusai2016,Phelan2018,Phelan2019,Phelan2020}. The Wiener-Hopf method is employed also to solve a large collection of mixed boundary value problems \citep{Duffy2008}. %\textcolor{red}{Check whether the citations Phelan 2018 should actually be Phelan 2019}

\section{Mathematical tools}

\subsection{Fourier transform and projection operators}\label{sec:Four_trans}

We define the Fourier transform of a function $f(x)$ as
\begin{equation}
\widehat{f}(\xi) = \mathcal{F}_{x\to\xi}[f(x)](\xi) = \int_{-\infty}^{+\infty} f(x) e^{i\xi x} dx\label{eq:forFour}
\end{equation}
and its inverse as
\begin{equation}
f(x) = \mathcal{F}_{\xi\to x}^{-1}[\widehat{f}(\xi)](x) = \frac{1}{2\pi} \int_{-\infty}^{+\infty} \widehat{f}(\xi) e^{-ix\xi} a\xi\label{eq:revFour}
\end{equation}
where $i$ is the imaginary unit. We are aware that it is advised to define Fourier space in terms of frequency $\nu$ \citep[see e.g][]{Press2007} rather than angular frequency or pulsation $\xi = 2\pi\nu$ (this is physics terminology if $x$ is interpreted as time). With $\nu$ the inverse transform lacks the factor $1/(2\pi)$, making it more symmetric with respect to the forward transform, and the Nyquist relation between grids in the normal and Fourier spaces simplifies to $\Delta x \Delta \nu = 1/N$, where $N$ is the number of grid points, without a factor $2\pi$ on the right-hand side. We also know that putting the minus sign in the exponent of the forward Fourier transform stresses the relationship with the Laplace transform and is the more common choice in fast Fourier transform (FFT) libraries, including the FFTW \citep{Frigo2005} used in \textsc{Matlab}. However, we chose the definition of the Fourier transform used normally in the two main application fields of Eq.~(\ref{eq:WHF}), physics and finance. An inconvenience is that the Fourier transform $\mathcal{F}$ of our equations translates into \texttt{ifft} in the \textsc{Matlab} code examples given in the supplementary material, and the inverse transform $\mathcal{F}^{-1}$ into \texttt{fft}.

We define the projection of a function $f(x)$ on the positive or on the negative half-axis through the multiplication with the indicator function of that set,
\begin{eqnarray}
\label{eq:P+}
f_+(x) &=& \mathcal{P}_{+,x}f(x) = 1_{\mathbb{R}_+}(x)f(x) \\
\label{eq:P-}
f_-(x) &=& \mathcal{P}_{-,x}f(x) = 1_{\mathbb{R}_-}(x)f(x).
\end{eqnarray}
A function that, like $f_+(x)$, is 0 for $x < 0$ and nonzero for $x > 0$ is called ``causal'', because it can be used to describe the effect of something that happens at $x = 0$ and causes the function to become nonzero. The two half-range Fourier transforms are
\begin{eqnarray}
\label{eq:FT+}
\widehat{f_+}(\xi) &=& \mathcal{F}_{x\to\xi}[f_+(x)](\xi) = \int_0^{+\infty} f(x) e^{i\xi x} dx \\
\label{eq:FT-}
\widehat{f_-}(\xi) &=& \mathcal{F}_{x\to\xi}[f_-(x)](\xi) = \int_{-\infty}^0 \, f(x) e^{i\xi x} dx.
\end{eqnarray}
The order in which the operators are applied matters:
\begin{equation}
\widehat{f_+}(\xi) = \mathcal{F}_{x\to\xi}[\mathcal{P}_{+,x}f(x)](\xi) \neq \widehat{f}_+(\xi) = \mathcal{P}_{+,\xi}\mathcal{F}_{x\to\xi}[f(x)](\xi).
\end{equation}
In other words, $\widehat{f_+}(\xi)$ is the Fourier transform of a function $f(x)$ that vanishes for negative arguments $x$, but $\widehat{f_+}(\xi)$ does not vanish itself for negative arguments $\xi$, which instead happens with $\widehat{f}_+(\xi)$; similarly for the $-$ case. The function $\widehat{f_+}(\xi)$ is analytic (or holomorphic), i.e.\ locally given by a convergent power series, in an upper complex half-plane that includes the real line; the function $\widehat{f_-}(\xi)$ is analytic in a lower complex half-plane that includes the real line. The half-range Fourier transforms can be considered special cases of the Laplace transform,
\begin{equation}
\widetilde{f}(s) = \mathcal{L}_{x\to s}[f(x)](s) = \int_0^{+\infty} f(x) e^{-sx} dx, \quad s \in \mathbb{C},
\end{equation}
where $s = \pm i\xi$, while the Fourier transform can be considered a special case of the bilateral or two-sided Laplace transform. Except possibly for $x = 0$, the indicator function $1_{\mathbb{R}_+}(x)$ coincides with the Heaviside step function $\Theta(x)$, and $1_{\mathbb{R}_-}(x)$ with $1-\Theta(x)$; $\Theta(x)$ is 1 if $x > 0$ and 0 if $x < 1$, while for $x = 0$ it can be assigned the value 0 (left-continuous choice), 1 (right-continuous choice), or 1/2 (symmetric choice). When integrating as in Eqs.~(\ref{eq:FT+}) and (\ref{eq:FT-}), the value for $x = 0$ matters only numerically and only if $x = 0$ is a grid point, as analytically the measure of a point is zero. Clearly the sum of the two projections, Eqs.~(\ref{eq:P+}) and (\ref{eq:P-}), is the full function,
\begin{equation}
f_+(x) + f_-(x) = f(x),
\end{equation}
and the sum of the two half-range Fourier transforms, Eqs.~(\ref{eq:FT+}) and (\ref{eq:FT-}), is the full Fourier transform,
\begin{equation}
\label{eq:PlemeljSokhotsky1a}
\widehat{f_+}(\xi) + \widehat{f_-}(\xi) = \widehat{f}(\xi).
\end{equation}

\subsection{Gibbs phenomenon}\label{sec:Gibbs_phenom}
As explained in Section \ref{sec:Four_trans}, we numerically implement the forward and inverse Fourier transform using the FFTW library in \textsc{Matlab}. The ranges of $x$ and $\xi$ cease to be infinite and continuous, and are approximated with grids of size $N$. The other parameter which defines both grids, which we centre around zero, is the truncation in the $x$ domain $x_{\max}$. The step is $\Delta x=2x_{\max}/N$ and the $x$ grid is
\begin{equation}\label{eq:xGrid}
x_n=n\Delta x,\quad n=-\frac{N}{2},-\frac{N}{2}+1,\dots,\frac{N}{2}-1.
\end{equation}
The step of the $\xi$ grid is given by the Nyquist relation, $\Delta\xi = 2\pi/(N\Delta x) = \pi/x_{\max}$; the truncation in the $\xi$ domain is $\xi_{\max}=\pi/\Delta x$ and the $\xi$ grid is
\begin{equation}\label{eq:xiGrid}
\xi_m=m\Delta\xi, \quad m=-\frac{N}{2},-\frac{N}{2}+1,\dots,\frac{N}{2}-1.
\end{equation}
The discrete forward and inverse Fourier transforms are
\begin{align}
\widehat{f}(\xi_m,\Delta x,N)&=\Delta x \sum^{N/2-1}_{n=-N/2}f\left(x_n\right)e^{i\xi_m x_n} \label{eq:DFT1}\\
f(x_n,\Delta\xi,N)&=\frac{\Delta \xi}{2\pi}\sum^{N/2-1}_{m=-N/2}\widehat{f}\left(\xi_m\right)e^{-ix_n\xi_m}.\label{eq:revDFT1}
\end{align}
The truncation of the sums in Eqs.~(\ref{eq:DFT1}) and (\ref{eq:revDFT1}) causes the Gibbs phenomenon. For a detailed explanation of its effect on the solution to Wiener-Hopf type equations see \cite{Phelan2019}. In this case we must consider two main issues: firstly, if the function $f(x)$ has a discontinuity, the truncation of $\widehat{f}(\xi_m,\Delta x,N)$ causes oscillations in $f(x_n,\Delta\xi,N)$ close to the discontinuity; secondly, the error away from that discontinuity will decay with the grid size $N$ as $|f(x_n)-f(x_n,\Delta\xi,N)| = O(1/N)$.

There have been many different approaches to solve or mitigate the Gibbs phenomenon; see e.g.\ \cite{vandeven1991family,gottlieb1997gibbs,tadmor2005adaptive,tadmor2007filters,ruijter2015application}. As in \cite{Phelan2019}, we apply a spectral filter in the Fourier domain, specifically the exponential filter of \cite{gottlieb1997gibbs}
\begin{align}
\label{eq:expFilt}
& \sigma(\eta)=e^{-\vartheta\eta^p},
\end{align}
where $p\in\mathbb{N}$ is even and $\eta = \xi/\xi_{\max}$. This function does not strictly meet the usual filter requirements described e.g.\ by \cite{vandeven1991family} as it does not go exactly to zero when $|\eta|=1$, nor do so its derivatives. However, if we select $\vartheta>-\log\varepsilon$, where $\varepsilon$ is machine precision, then the filter coefficients are within computational accuracy of the requirements. Advantages of the exponential filter are its good performance, its simple form, and the order of the filter being equal to the parameter $p$ which is directly input into the filter equation.

We also investigated the use of the Planck taper \citep{mckechan2010tapering}, which is defined piecewise as
\begin{align}
&\sigma(\eta)=
\begin{cases}
0, &\eta\leq \eta_1,\ \ \qquad \eta_1=-1 \\
\frac{1}{e^{\;z(\eta)}+1},\ z(\eta)=\frac{\eta_2-\eta_1}{\eta-\eta_1}+\frac{\eta_2-\eta_1}{\eta-\eta_2}, & \eta_1<\eta<\eta_2,\ \, \eta_2=\epsilon-1\\
1, & \eta_2 \leq\eta\leq\eta_3,\ \,\eta_3=1-\epsilon\\
\frac{1}{e^{\;z(\eta)}+1},\ z(\eta)=\frac{\eta_3-\eta_4}{\eta-\eta_3}+\frac{\eta_3-\eta_4}{\eta-\eta_4}, &\eta_3<\eta<\eta_4,\ \, \eta_4=1\\
0, &\eta\geq\eta_4.
\end{cases}
\end{align}
Here, the value of $\epsilon$ gives the proportion of the range of $\eta$ which is used for the slope regions.
In common with the findings by \cite{Phelan2019}, the Planck taper, whilst having some interesting characteristics such as a flat central section and a filter order of $\infty$, when tested did not offer any advantage over the exponential filter, so we did not pursue its use any further.

\subsection{Hilbert transform and Wiener-Hopf factorisation}

The Hilbert transform \citep[see e.g][]{Pandey1996,VergaraCaffarelli1999,King2009} of $f(x)$ is defined as the Cauchy principal value of the convolution with $1/(\pi x)$,
\begin{eqnarray}
\label{eq:Hilbert_transform}
\mathcal{H}_x f(x) &=& \mathrm{p.v.}\, \frac{1}{\pi x} * f(x) = \mathrm{p.v.}\, \frac{1}{\pi} \int_{-\infty}^{+\infty} \frac{f(x')}{x-x'} dx' \\
&=& \lim_{\epsilon \to 0^+} \frac{1}{\pi} \left( \int_{-1/\epsilon}^{-\epsilon} \frac{f(x')}{x-x'} dx' + \int_\epsilon^{1/\epsilon} \frac{f(x')}{x-x'} dx' \right).
\end{eqnarray}
The principal value avoids that the improper integral evaluates to the indefinite form $+\infty-\infty$. Since with this definition the Hilbert transform often appears multiplied by the imaginary unit (see the following equations), some authors such as \cite{Stenger1973} define the Hilbert transform as the principal value of the convolution with $i/(\pi x)$. To stress that the Hilbert transform is actually a functional like the Fourier and Laplace transforms, one could write $\mathcal{H}_{x'\to x}[f(x')](x)$ instead than $\mathcal{H}_x f(x)$, but this is too cumbersome and also less useful because, unlike $x$ and $\xi$ or $x$ and $s$, $x'$ and $x$ belong to the same space. A subscript like $x \to \xi$ or $x$ can be omitted when there is no misunderstanding about which variable the operators $\mathcal{F, F}^{-1}, \mathcal{P_+, P_-}$ and $\mathcal{H}$ act on, notably when the argument function depends on a single variable; this is mostly the case here. The operator $i\mathcal{H}$ is its own inverse,
\begin{equation}
\label{eq:HT_inversion}
(i\mathcal{H})^2 f(x) = f(x);
\end{equation}
equivalently, $\mathcal{H}^{-1} = -\mathcal{H}$. The convolution theorem
\begin{equation}
\label{eq:convolution_theorem}
\mathcal{F}[f(x)*g(x)] = \widehat{f}(\xi)\widehat{g}(\xi),
\end{equation}
which maps the convolution product in real space to a simple product in Fourier space, together with the Fourier transform \citep{Weisstein2021}
\begin{equation}
\label{eq:FT_inverse_function}
\mathrm{p.v.}\,\mathcal{F}\frac{1}{\pi x} = i\,\mathrm{sgn}\,\xi
\end{equation}
enables to express the Hilbert transform through a forward and an inverse Fourier transform,
\begin{equation}
\label{eq:HilbertFourier}
\mathcal{H}f(x) = \mathcal{F}^{-1}\big[i\,\mathrm{sgn}\,\xi\,\widehat{f}(\xi)\big].
\end{equation}
Thus a fast method to compute the Hilbert transform numerically consists simply in evaluating Eq.~(\ref{eq:HilbertFourier}) through a forward and an inverse FFT. In the next subsection, \ref{sec:HilbertSinc}, we shall see more sophisticated numerical methods. Swapping the forward Fourier transform with its inverse, Eq.~(\ref{eq:FT_inverse_function}) becomes
\begin{equation}
\mathrm{p.v.}\,\mathcal{F}^{-1}\frac{1}{\pi\xi} = -i\,\mathrm{sgn}\,x,
\end{equation}
and Eq.~(\ref{eq:HilbertFourier}) becomes
\begin{equation}
\label{eq:HilbertFourier2}
i\mathcal{H}\widehat{f}(\xi) = \mathcal{F}\big[\mathrm{sgn}\,x\,f(x)\big];
\end{equation}
substituting $\mathrm{sgn}\,x = 1_{\mathbb{R}_+}(x) - 1_{\mathbb{R}_-}(x)$ (this is true also for $x = 0$, while $\mathrm{sgn}\,x = 2\Theta(x) - 1$ is fulfilled for $x = 0$ only if $\Theta(0) = 1/2$), and applying the definitions of the half-range Fourier transforms, Eqs.~(\ref{eq:FT+}) and (\ref{eq:FT-}), yields
\begin{equation}
\label{eq:PlemeljSokhotsky1b}
\widehat{f_+}(\xi) - \widehat{f_-}(\xi) = i\mathcal{H}\widehat{f}(\xi).
\end{equation}
This can be shown also evaluating the integral in Eq.~(\ref{eq:Hilbert_transform}) with contour integration methods in the complex plane. Together, Eqs.~(\ref{eq:PlemeljSokhotsky1a}) and (\ref{eq:PlemeljSokhotsky1b}) are known as Plemelj-Sokhotsky relations \citep{Pandey1996,VergaraCaffarelli1999,King2009}. They can be rearranged as
\begin{eqnarray}
\label{eq:PlemeljSokhotsky2a}
\widehat{f_+}(\xi) &=& \frac{1}{2}\big[\widehat{f}(\xi) + i\mathcal{H}\widehat{f}(\xi)\big] \\
\label{eq:PlemeljSokhotsky2b}
\widehat{f_-}(\xi) &=& \frac{1}{2}\big[\widehat{f}(\xi) - i\mathcal{H}\widehat{f}(\xi)\big]
\end{eqnarray}
or, with a different notation involving the Fourier-transform and projection operators, as
\begin{eqnarray}
\label{eq:PlemeljSokhotsky3a}
\mathcal{F}\mathcal{P_+}f(x) &=& \frac{1}{2}\big[\mathcal{F}f(x) + i\mathcal{H}\mathcal{F}f(x)\big] \\
\label{eq:PlemeljSokhotsky3b}
\mathcal{F}\mathcal{P_-}f(x) &=& \frac{1}{2}\big[\mathcal{F}f(x) - i\mathcal{H}\mathcal{F}f(x)\big].
\end{eqnarray}
%\textit{From here up to the beginning of Eq. 2.35 $f$ and $\widehat{f}$ did not have $(x)$ and $(\xi)$ written explicitly which is inconsistent with the other equations so I have added these in.}
Substituting $f(x)$ with $\mathcal{P}_+f(x)$ in Eq.~(\ref{eq:PlemeljSokhotsky3a}) and $f(x)$ with $\mathcal{P}_-f(x)$ in Eq.~(\ref{eq:PlemeljSokhotsky3b}), and taking into account that projection operators are idempotent, i.e., $\mathcal{PP}f (x)= \mathcal{P}f(x)$, shows that the half-range Fourier transforms are eigenfunctions of the Hilbert transform operator,
\begin{eqnarray}
i\mathcal{H}\widehat{f_+}(\xi) & = & \phantom{-}\widehat{f_+}(\xi) \label{eq:Heigenf+} \\
i\mathcal{H}\widehat{f_-} (\xi)& = & -\widehat{f_-}(\xi). \label{eq:Heigenf-}
\end{eqnarray}
This is evident also substituting $f(x)$ with $f_+(x)$ or with $f_-(x)$ in Eq.~(\ref{eq:HilbertFourier2}), or applying the operator $i\mathcal{H}$ to both sides of Eqs.~(\ref{eq:PlemeljSokhotsky2a})--(\ref{eq:PlemeljSokhotsky2b}) and simplifying with Eq.~(\ref{eq:HT_inversion}). Eqs.~(\ref{eq:Heigenf+}) and (\ref{eq:Heigenf-}) allow to obtain Eq.~(\ref{eq:PlemeljSokhotsky1b}) applying the operator $i\mathcal{H}$ to both sides of Eq.~(\ref{eq:PlemeljSokhotsky1a}); conversely, Eq.~(\ref{eq:PlemeljSokhotsky1a}) can be reobtained applying $i\mathcal{H}$ to both sides of Eq.~(\ref{eq:PlemeljSokhotsky1b}). Eqs.~(\ref{eq:PlemeljSokhotsky2a}) and (\ref{eq:PlemeljSokhotsky2b}) are invariant with respect to an application of $i\mathcal{H}$ to both sides.

The key step in the Wiener-Hopf solution of Eq.~(\ref{eq:WHF}) described in the following section is the decomposition of a function $\widehat{f}$, i.e., the reverse of Eq.~(\ref{eq:PlemeljSokhotsky1a}),
\begin{equation}
\widehat{f}(\xi) = \widehat{f_+} (\xi)+ \widehat{f_-}(\xi).
\end{equation}
The factorisation of a function $\widehat{f}(\xi)$
\begin{equation}
\widehat{f}(\xi) = \widehat{f_+}(\xi) \widehat{f_-}(\xi),
\end{equation}
which is required too, can be reduced to a decomposition by taking logarithms,
\begin{equation}
\log\widehat{f}(\xi) = \log\widehat{f_+}(\xi) + \log\widehat{f_-}(\xi);
\end{equation}
this procedure is called logarithmic decomposition. As described by \cite{Rino1970}, \cite{Henery1974} and \cite{Bart2004}, the decomposition can be achieved by
\begin{eqnarray}
\label{eq:Rino_decomposition_+}
\widehat{f_+}(\xi) &=& \mathcal{FP_+F}^{-1}\widehat{f}(\xi) \\
\label{eq:Rino_decomposition_-}
\widehat{f_-}(\xi) &=& \mathcal{FP_-F}^{-1}\widehat{f}(\xi),
\end{eqnarray}
as can also be seen from the definitions of the half-range Fourier transforms, Eqs.~(\ref{eq:FT+}) and (\ref{eq:FT-}). More generally by the Plemelj-Sokhotsky relations, Eqs.~(\ref{eq:PlemeljSokhotsky2a}) and (\ref{eq:PlemeljSokhotsky2b}) can be used \citep{Stenger1973}. Eqs.~(\ref{eq:Rino_decomposition_+}) and (\ref{eq:Rino_decomposition_-}) are a special case of the latter if the Hilbert transform is computed through Eq.~(\ref{eq:HilbertFourier2}). % the two approaches yield identical results.

The definition of the two half-range Fourier transforms, Eqs.~(\ref{eq:FT+}) and (\ref{eq:FT-}), can be generalized splitting the $x$ axis around a constant $a \neq 0$. \cite{Feng2008} showed how the shift theorem,
\begin{equation}
\mathcal{F}_{x\to\xi}[f(x+a)](\xi) = \widehat{f}(\xi)e^{-ia\xi},
\end{equation}
 can be exploited to generalise the Plemelj-Sokhotsky relations to
\begin{eqnarray}
\label{eq:PlemeljSokhotsky4a}
\widehat{f_+}(\xi) &=& % e^{ia\xi}\big(\widehat{f}(\xi)e^{-ia\xi}\big)_+ =
\frac{1}{2}\big\{\widehat{f}(\xi) + e^{ia\xi}i\mathcal{H}_\xi\big[\widehat{f}(\xi)e^{-ia\xi}\big]\big\} \\
\label{eq:PlemeljSokhotsky4b}
\widehat{f_-}(\xi) &=& % e^{ia\xi}\big(\widehat{f}(\xi)e^{-ia\xi}\big)_- =
\frac{1}{2}\big\{\widehat{f}(\xi) - e^{ia\xi}i\mathcal{H}_\xi\big[\widehat{f}(\xi)e^{-ia\xi}\big]\big\}.
\end{eqnarray}
It might be a good idea to write $\widehat{f_{+,a}}$ and $\widehat{f_{-,a}}$ on the left-hand side, but we will not do it to not overburden the notation, since it will be clear from the context with respect to which parameter a function is decomposed. In the above formulas Eq.~(\ref{eq:HilbertFourier2}) generalizes to
\begin{eqnarray}
e^{ia\xi}i\mathcal{H}_\xi\big[\widehat{f}(\xi)e^{-ia\xi}\big] & = & e^{ia\xi}\mathcal{F}_{x\to\xi}\big[\mathrm{sgn}(x)f(x+a)\big] \\
& = & \mathcal{F}_{x\to\xi}\big[\mathrm{sgn}(x-a)f(x)\big] \\
& = & \mathcal{F}_{x\to\xi}\big[(1_{(a,+\infty)}(x) - 1_{(-\infty,a)}(x))f(x)\big],
\end{eqnarray}
and thus it is easy to show that
\begin{gather}
\lim_{a\to-\infty} e^{ia\xi}i\mathcal{H}_\xi\big[\widehat{f}(\xi)e^{-ia\xi}\big] = \widehat{f}(\xi) \\
\lim_{a\to+\infty} e^{ia\xi}i\mathcal{H}_\xi\big[\widehat{f}(\xi)e^{-ia\xi}\big] = -\widehat{f}(\xi)
\end{gather}
and that $\lim_{a\to-\infty} \widehat{f_+}(\xi) = \widehat{f}(\xi),\ \lim_{a\to+\infty} \widehat{f_+}(\xi) = 0,\ \lim_{a\to-\infty} \widehat{f_-}(\xi) = 0,\ \lim_{a\to+\infty} \widehat{f_-}(\xi) = \widehat{f}(\xi)$. These limits are useful to retrieve the results for the classical Wiener-Hopf equation from those for the Fredholm equation.

\subsection{Fast Hilbert transform with sinc functions}
\label{sec:HilbertSinc}
%Eq.~(\ref{eq:HilbertFourier}) provides a straightforward method to evaluate numerically the Hilbert transform. Since it is based on two FFTs, its computational cost is $O(N\log N)$, where $N$ is the number of grid points, and thus is called fast. It yields exactly the same results as the following $O(N^2)$ quadrature method \cite{King2009}, and in both cases the error converges polynomially with n, so clearly, among these two, the FFT approach is preferable because of its speed. The quadrature method consists in {\textcolor{red}...} For details, see the function htq.m in the supplementary material.

Eq.~(\ref{eq:HilbertFourier}) provides a straightforward method to evaluate numerically the Hilbert transform. Since it is based on two FFTs, its computational cost is $O(N\log N)$, where $N$ is the number of grid points, and thus is called fast. We compared this method with the quadrature method described in \citet[Eqs.~(4.19)--(4.20)]{King2009}, where the summation is taken over every second point in order to avoid the singularity which results when $x_i-x_j=0$. We tried various quadrature weights, including Simpson's rule and 3rd and 4th order quadrature \citep{Press2007}. For our implementation, see the \textsc{Matlab} functions htq.m and weights.m in the supplementary material. All weights give the same result and have polynomial convergence with $N$. Therefore, as with quadrature the computation speed is $O(N^2)$, the FFT-based method is preferable because of its higher speed.

An alternative, but equally fast $O(N\log N)$ approach to compute numerically the Hilbert transform is based on the sinc expansion approximation of analytical functions. Sinc functions are deeply studied in the books by \cite{Stenger1993,Stenger2011}, where the author shows that a function $f$ analytical in the whole complex plane and of exponential type with parameter $\pi/h$, i.e.,
\begin{equation}
|f(z)|\leq C e^{\pi|z|/h},
\end{equation}
can be reconstructed exactly from the knowledge of its values on an equispaced grid of step $h$. Defining the sinc functions
\begin{equation}
S_n(z,h) = \frac{\sin(\pi(z-nh)/h)}{\pi(z-nh)/h}, \quad n \in \mathbb{Z},
\end{equation}
the function $f$ admits the sinc expansion \cite[Theorem 1.10.1]{Stenger1993}
\begin{equation}
f(z) = \sum_{n=-\infty}^{+\infty} f(nh) S_n(z,h).
\end{equation}
Moreover, also its Fourier transform admits the sinc expansion
\begin{equation}
\widehat{f}(\xi) = h\sum_{n=-\infty}^{+\infty} f(nh) e^{i\xi nh} \quad \textrm{if}\ |\xi|<\pi/h,
\end{equation}
while it is zero if $|\xi|\geq\pi/h$. Finally, also the integrals of $f$ and $|f|^2$ can be written as sinc expansions,
\begin{equation}
\int_{-\infty}^{+\infty} f(x) dx = h\sum_{n=-\infty}^{+\infty} f(nh), \quad \int_{-\infty}^{+\infty} |f(x)|^2 dx = h\sum_{n=-\infty}^{+\infty} |f(nh)|^2.
\end{equation}
The above results show in particular that the trapezoidal quadrature rule with step size $h$ is exact. Using the following result on the Hilbert transform of the sinc functions \cite[Corollary 6.1]{Feng2008}
\begin{equation}
\mathcal{H}_\xi S_n(\xi,h) = \frac{1-\cos(\pi(\xi-nh)/h)}{\pi(\xi-nh)/h},
\end{equation}
also the Hilbert transform can be evaluated exactly as
\begin{equation}
\mathcal{H}_\xi f(\xi) = \sum_{n=-\infty}^{+\infty} f(nh)\frac{1-\cos(\pi(\xi-nh)/h)}{\pi(\xi-nh)/h}.\label{eq:Hilbfinal}
\end{equation}
This holds for a function $f$ that is analytic in the whole complex plane. However, this can be used also to approximate a function that is analytic only in a strip including the real axis, which is the case considered in this article. As shown in \citet[Theorems 3.1.3, 3.2.1 and 3.1.4]{Stenger1993}, in this case the trapezoidal approximation has a discretisation error that decays exponentially with respect to $h$. \citet[Section 6.5]{Feng2008} also shows that the computation of the Hilbert transform via a sinc expansion can be performed using the FFT. The idea is that to compute a discrete Hilbert transform it is necessary to do matrix-vector multiplications involving Toeplitz matrices. As is well known, this kind of multiplications can be performed exploiting the FFT, once those matrices are embedded in a circulant matrix as in \citet[Appendix B]{Feng2008} and \cite{Fusai2012}. %In particular, Feng and Linetsky, with regard to the computation of the Hilbert transform \citep[Theorem 3.3]{Feng2009} and of the whole Plemelj-Sokhotski equations (\ref{eq:PlemeljSokhotsky2a})--(\ref{eq:PlemeljSokhotsky2b}) \cite[Theorem 6.5]{Feng2008} \cite[Theorem 3.4]{Feng2009} with sinc functions, prove the following convergence result: if a function is analytic in a suitable strip around the real axis, then the discretisation error of its numerical factorisation or decomposition decays exponentially with respect to the number of discretisation points $N$. Further details can be found in the cited references \cite[Section 6]{Feng2008} \cite[Section 3.4]{Feng2009}. An implementation of both the $O(N\log N)$ methods presented above can be found in the fuction ifht.m in the supplementary material.
In particular, \cite{Feng2008,Feng2009} prove the following convergence result: if a function is analytic in a suitable strip around the real axis, then the discretisation error of its numerical factorisation or decomposition decays exponentially with respect to the step size between discretisation points $h$. \citet[Theorem 3.3]{Feng2009} considers the computation of the Hilbert transform and \citet[Theorems 6.5]{Feng2008} and \citet[Theorems 3.4]{Feng2009} are concerned with the calculation of the Plemelj-Sokhotsky relations from Eqs.~(\ref{eq:PlemeljSokhotsky2a}) and (\ref{eq:PlemeljSokhotsky2b}) in particular. An implementation of both the $O(N\log N)$ methods presented above can be found in the fuction ifht.m in the supplementary material.

However, %as also explained by \citet[Theorem ?? ]{Stenger1993},
the error due to the infinite sum in Eq.~(\ref{eq:Hilbfinal}) being truncated to the number of FFT grid points depends on the shape of the function under transform. This has been explored further in the error bounds developped in \citet[Section 6.4.2]{Feng2008} and \citet[Section 3]{Phelan2019}. In the case of a function that decays exponentially as $|\xi|\rightarrow \infty$, exponential convergence is obtained. However, if a function decays polynomially then the truncation error convergence is only polynomially decreasing. In their paper on lookback options \cite[Theorem 3.3]{Feng2009} report a result by Stenger that proves the exponential convergence of the discrete sinc-based fast Hilbert transform to the continuous Hilbert transform. They also examine the truncation error, specifically observing in a footnote that this converges exponentially only under certain conditions, notably $f(x)\leq\kappa \exp(c|x|^\nu)$ for some $\kappa, c, \nu>0$. This algorithm can be obtained with an eigenfunction expansion of $\mathcal{H}$ and is identical to the Kress and Martensen method, which was introduced with a proof that its error converges exponentially \citep{Kress1970,Weideman1995}.

%\textcolor{red}{\textit{See also Feng and Linetsky's paper on lookback options \cite[Theorem 10]{Feng2009} {\color{blue}(CEP note: which theorem, the published paper is labelled differently)}, where they report a result by Stenger that proves the exponential convergence of the discrete Hilbert transform with sinc functions to the continuous Hilbert transform {\color{blue} although the truncation error only converges exponentially under certain conditions, notably that $f(x)\leq\kappa \exp(c|x|^\nu)$ for some $\kappa, c, \nu>0$}. Mention that this algorithm can be obtained with an eigenfunction expansion of $\mathcal{H}$ and is identical to the Kress and Martensen method, which was introduced with a proof that its error converges exponentially \cite{Kress1970,Weideman1995}.}}

\section{The classical Wiener-Hopf method for a convolution equation on a semi-infinite interval}
\label{sec:CWHE}

Consider Eq.~(\ref{eq:WHF}) with $b = +\infty$. The lower integration limit $a$ can be set to 0 shifting the $x$ scale horizontally by the constant $a$ to a new scale $x' = x - a$; the prime is dropped hereafter. The functions $f(x)$ and $g(x)$, whose domain is $[0,+\infty)$, are extended to the whole real axis defining $f_0(x) = 0$ for $x < 0$, $f_0(x) = f(x)$ for $ \ge 0$ and $g_0(x) = 0$ for $x < 0$, $g_0(x) = g(x)$ for $x \ge 0$. Define moreover the auxiliary function
\begin{equation}
f_1(x) = \int_0^{\infty} k(x-x') f(x') dx' = \int_{-\infty}^{+\infty} k(x-x') f_0(x') dx', \quad x < 0,
\end{equation}
and $f_1(x) = 0$ for $x \ge 0$, i.e., $\widehat{f_1} = \big(\widehat{k}\widehat{f_0}\big)_-$. Since $f_0$ and $g_0$ are $+$ functions and $f_1$ is a $-$ function, % i.e., $(f_0)_+ = f_0$, $(g_0)_+ = g_0$, and $(f_1)_- = f_1$,
it is customary to denote these functions $f_+$, $g_+$, and $f_-$ respectively. With them Eq.~(\ref{eq:WHF}) is extended to
\begin{equation}
\label{eq:WH1_extended}
\lambda f_+(x) - \int_{-\infty}^{+\infty} k(x-x') f_+(x') dx' + f_-(x) = g_+(x),
\quad x \in \mathbb{R},
\end{equation}
or, with a more compact notation for the convolution,
\begin{equation}
\lambda f_+(x) - k(x)*f_+(x) + f_-(x) = g_+(x), \quad x \in \mathbb{R}.
\end{equation}
The extension of the integration domain to the whole real axis does not affect the equation and its solution on the positive half-axis, but allows to apply the convolution theorem, Eq.~(\ref{eq:convolution_theorem}), to obtain the equation in Fourier space,
\begin{equation}
\label{eq:WH1_FT}
\widehat{l}(\xi)\widehat{f_+}(\xi) + \widehat{f_-}(\xi) = \widehat{g_+}(\xi),
\end{equation}
where $\widehat{l}(\xi) = \lambda-\widehat{k}(\xi)$, ($l(x)$ and $\widehat{l}(\xi)$ are the functional derivatives of the equation with respect to the solution in normal and Fourier space respectively). Dropping the argument $\xi$ for brevity, factorising $\widehat{l} = \widehat{l_-}\widehat{l_+}$ and dividing the equation by $\widehat{l_-}$ gives
\begin{equation}
\label{eq:WH1_FT_2}
\widehat{l_+}\widehat{f_+} + \widehat{l_-}^{-1}\widehat{f_-} = \widehat{l_-}^{-1}\widehat{g_+}.
\end{equation}
Defining
\begin{equation}
\widehat{c} = \widehat{l_-}^{-1}\widehat{g_+}
\end{equation}
and decomposing it as $\widehat{c} = \widehat{c_+}+ \widehat{c_-}$ yields finally
\begin{equation}
\widehat{l_+}\widehat{f_+} + \widehat{l_-}^{-1}\widehat{f_-}= \widehat{c_+}+ \widehat{c_-}.
\end{equation}
The $+$ and $-$ components can be separated into
\begin{gather}
\widehat{f_+} = \widehat{l_+}^{-1}\widehat{c_+} \\
\widehat{f_-} = \widehat{l_-}\widehat{c_-},
\end{gather}
which allows to obtain the sought solution from
\begin{equation}
\label{eq:WH1_solution}
f_+ = \mathcal{F}^{-1}\big( \widehat{l_+}^{-1}\widehat{c_+} \big),
\end{equation}
while $f_-$ was introduced as an auxiliary function and is not of further interest.

The case with $a = -\infty$ is treated in a similar fashion. The upper integration limit $b$ can be set to 0 shifting the $x$ scale horizontally by the constant $b$ to a new scale $x' = x - b$; the prime is dropped hereafter. The functions $f(x)$ and $g(x)$, whose domain is $(-\infty,0]$, are extended to the whole real axis defining $f_0(x) = f(x)$ for $x \le 0$, $f_0(x) = 0$ for $x > 0$ and $g_0(x) = g(x)$ for $x \le 0$, $g_0(x) = 0$ for $x > 0$. Define moreover the auxiliary function
\begin{equation}
f_2(x) = \int_{-\infty}^0 k(x-x') f(x') dx' = \int_{-\infty}^{+\infty} k(x-x') f_0(x') dx', \quad x > 0,
\end{equation}
and $f_2(x) = 0$ for $x \le 0$, i.e., $\widehat{f_2} = \big(\widehat{k}\widehat{f_0}\big)_+$. Now $f_0$ and $g_0$ are $-$ functions and $f_2$ is a $+$ function, % i.e., $(f_0)_- = f_0$, $(g_0)_- = g_0$, and $(f_2)_+ = f_2$,
so it is customary to denote these functions $f_-$, $g_-$, and $f_+$ respectively. With them Eq.~(\ref{eq:WHF}) is extended to equations identical to Eqs.~(\ref{eq:WH1_extended})--(\ref{eq:WH1_solution}), except that the $+$ and $-$ indices are swapped. In particular, the sought solution is obtained from
\begin{gather}
\widehat{c} = \widehat{l_+}^{-1}\widehat{g_-} \\
f_- = \mathcal{F}^{-1}\big( \widehat{l_-}^{-1} \widehat{c_-} \big).
\end{gather}

A more elegant alternative to shifting the $x$ scale forth and back by the constant $a$ or $b$ is to modulate the functions in Fourier space decomposing $\widehat{c}$ with respect to this constant by the generalized Plemelj-Sokhotsky relations, Eqs.~(\ref{eq:PlemeljSokhotsky4a}) and (\ref{eq:PlemeljSokhotsky4b}). The function $\widehat{l}$ is always factorized with respect to 0, while $\widehat{c}$ is decomposed with respect to $a$ when $b = +\infty$ or to $b$ when $a = -\infty$. For details, see the function whf\_gmf\_filt4.m in the supplementary material. % Sec.~\ref{sec:whf_gmf.m}.

% the method is applicable also to an asymmetric and complex kernel, as long as its integral is smaller than one

\section{Generalisation of the Wiener-Hopf method to a convolution equation on a finite interval: the Fredholm equation}
\label{sec:FE}

\subsection{Theory}

In the Fredholm equation both integration limits $a$ and $b$ are finite; either $a$ or, less usually, $b$ can be set to 0 shifting the $x$ scale, but, unlike with the classical Wiener-Hopf equation described in the previous section, we prefer not to modify any of the two integration limits; instead, we will use the generalized Plemelj-Sokhotsky relations. The functions $f(x)$ and $g(x)$, whose domain is $[a,b]$, are extended to the whole real axis defining $f_0(x) = f(x)$ for $x \in [a,b]$, $f_0(x) = 0$ for $x \notin [a,b]$ and $g_0(x) = g(x)$ for $x \in [a,b]$, $g_0(x) = 0$ for $x \notin [a,b]$. The kernel $k(x)$, whose domain is $[a-b,b-a]$, is extended to the whole real axis defining $k_0(x) = k(x)$ for $x \in [a-b,b-a]$ and $k_0(x) = 0$ for $x \notin [a-b,b-a]$. Define moreover the two auxiliary functions
\begin{equation}
f_1(x) = \int_a^b k(x-x') f(x') dx' = \int_{-\infty}^{+\infty} k_0(x-x') f_0(x') dx', \quad x < a, \label{eq:F1_def}
\end{equation}
$f_1(x) = 0$ for $x \ge a$, i.e., $\widehat{f_1} = e^{ia\xi}\big(\widehat{k_0}\widehat{f_0}e^{-ia\xi}\big)_- = \widehat{f_-}$,
\begin{equation}
f_2(x) = \int_a^b k(x-x') f(x') dx' = \int_{-\infty}^{+\infty} k_0(x-x') f_0(x') dx', \quad x > b,\label{eq:F2_def}
\end{equation}
and $f_2(x) = 0$ for $x \le b$, i.e., $\widehat{f_2} = e^{ib\xi}\big(\widehat{k}\widehat{f_0}e^{-ib\xi}\big)_+ = \widehat{f_+}$. Because $k_0(x) = 0$ for $x \notin [a-b,b-a]$, $f_-(x) = 0$ also for $x < a-(b-a) = 2a-b$ and $f_+(x) = 0$ also for $x > b-(a-b) = 2b-a$. Thus Eq.~(\ref{eq:WHF}) extends to
\begin{equation}
\label{eq:F_extended}
\lambda f_0(x) - \int_{-\infty}^{+\infty} k_0(x-x') f_0(x') dx' + f_-(x) + f_+(x) = g_0(x)
\end{equation}
or, with a more compact notation for the convolution,
\begin{equation}
\lambda f_0(x) - k_0(x)*f_0(x) + f_-(x) + f_+(x) = g_0(x),
\end{equation}
and upon Fourier transformation, setting $\widehat{l}(\xi) = \lambda - \widehat{k_0}(\xi)$,
\begin{equation}
\label{eq:F_FT}
\widehat{l}(\xi)\widehat{f_0}(\xi) + \widehat{f_-}(\xi) + \widehat{f_+}(\xi) = \widehat{g_0}(\xi).
\end{equation}
Eqs.~(\ref{eq:F_extended})--(\ref{eq:F_FT}) look similar to Eqs.~(\ref{eq:WH1_extended})--(\ref{eq:WH1_FT}), but now we have two auxiliary functions, $\widehat{f_-}(\xi)$, which is $-$ with respect to any $c \ge a$, and $\widehat{f_+}(\xi)$, which is $+$ with respect to any $d \le b$, while both the unknown function $\widehat{f_0}(\xi)$ and the forcing function $\widehat{g_0}(\xi)$ are $+$ with respect to $a$ (or any number $\le a$) and $-$ with respect to $b$ (or any number $\ge b$). Therefore the usual approach is to split Eq.~(\ref{eq:F_FT}) into two coupled Wiener-Hopf equations, one with the origin shifted to $a$, the other with the origin shifted to $b$ \citep{Green2010}. These functions typically involve the four redundant unknowns $\widehat{f_0}(\xi)e^{-ia\xi}$, $\widehat{f_+}(\xi)e^{-ib\xi}$ (which are $+$ functions), $\widehat{f_0}(\xi)e^{-ib\xi}$, $\widehat{f_-}(\xi)e^{-ia\xi}$ (which are $-$ functions). In Section \ref{sec:iterative} the functions $\widehat{f_-}(\xi)$ and $\widehat{f_+}(\xi)$ correspond to $J_-$ and $J_+$ from \cite{Green2010}, Eq.~(2.51), while $\widehat{c}_1$ and $\widehat{c}_2$ correspond to $P$ and $Q$ from Eqs.~(2.12) and (2.24) in that paper.

\subsection{Iterative solution}
\label{sec:iterative}

We solved the system of integral equations described in Eqs.~(\ref{eq:F1_def})--(\ref{eq:F_extended}) iteratively observing that, if we know $\widehat{f_+}(\xi)$ and subtract it from both sides of Eq.~(\ref{eq:F_FT}) with the origin of the $x$ axis shifted to $a$, the result looks like Eq.~(\ref{eq:WH1_FT}), so that we can use the method described in Sec.~\ref{sec:CWHE} to obtain $\widehat{f_-}(\xi)$; similarly, if we know $\widehat{f_-}(\xi)$ and subtract it from both sides of Eq.~(\ref{eq:F_FT}) with the origin of the $x$ axis shifted to $b$, we can use the method described in Sec.~\ref{sec:CWHE} to obtain $\widehat{f_+}(\xi)$. Thus, once again dropping the argument $\xi$ for brevity of notation, our procedure is to write Eq.~(\ref{eq:F_FT}) divided once by $\widehat{l_-}$, as in Eq.~(\ref{eq:WH1_FT_2}), and once by $\widehat{l_+}$,
\begin{gather}
\widehat{l_+}\widehat{f_0} + \widehat{l_-}^{-1}\widehat{f_-} + \widehat{l_-}^{-1}\widehat{f_+} = \widehat{l_-}^{-1}\widehat{g_0} \label{eq:F_FT_i1} \\
\widehat{l_-}\widehat{f_0} + \widehat{l_+}^{-1}\widehat{f_-} + \widehat{l_+}^{-1}\widehat{f_+} = \widehat{l_+}^{-1}\widehat{g_0}, \label{eq:F_FT_i2}
\end{gather}
start from the guess $\widehat{f_+} = 0$ in Eq.~(\ref{eq:F_FT_i1}), set
\begin{equation}
\widehat{c_1} = \widehat{l_-}^{-1}(\widehat{g_0}-\widehat{f_+}),
\end{equation}
decompose $\widehat{c_1} = \widehat{c_{1+}} + \widehat{c_{1-}}$ with respect to $a$, and compute the approximations
\begin{gather}
\widehat{f_0} = \widehat{l_+}^{-1}\widehat{c_{1+}}\\
\widehat{f_-} = \widehat{l_-}\widehat{c_{1-}}
\end{gather}
as $+$ and $-$ functions with respect to $a$; then turn to Eq.~(\ref{eq:F_FT_i2}), set
\begin{equation}
\widehat{c_2} = \widehat{l_+}^{-1}(\widehat{g_0}-\widehat{f_-}),
\end{equation}
decompose $\widehat{c_2} = \widehat{c_{2+}} + \widehat{c_{2-}}$ with respect to $b$, and compute the new approximations
\begin{gather}
\widehat{f_0} = \widehat{l_-}^{-1}\widehat{c_{2-}}\\
\widehat{f_+} = \widehat{l_+}\widehat{c_{2+}}
\end{gather}
as $+$ and $-$ functions with respect to $b$; and so on until the difference between the values of $\widehat{f_0}$ at an iteration and the previous falls below a threshold. An equivalent result is obtained starting from the guess $\widehat{f_-} = 0$ in Eq.~(\ref{eq:F_FT_i2}) and the computation of $\widehat{c_2}$. Notice that the iterations are performed looking for a fixed point on the variables $\widehat{f_-}$ and $\widehat{f_+}$, while $\widehat{f_0}$ is a side product output at each step, but not used to compute the next step. For details, see the function whf\_gmf\_filt4.m in the supplementary material. % Sec.~\ref{sec:whf_gmf.m}.

%\subsection{Henery's iterative solution}\label{sec:iterative_Henery}
\subsection{Other iterative solutions}\label{sec:iterative_other}
In this journal, \cite{Henery1977} proposed an iterative solution of the Fredholm equation, but presented only the theory without a numerical validation. In our tests, its literal implementation does not work. The procedure can be mapped to ours including a missing projection and an untold detail: the $y_n$ found in the residual correction scheme are corrections to the solution and thus must be added together. %\textcolor{red}{\textit{Details in henery.tex.}
However, \cite{Henery1977} did not express the algorithm in terms of the Hilbert transform and thus used only the simple implementation with the sign function, but not with a sinc function expansion.

\cite{Margetis2006} presented an iterative solution limited to algebraic kernel functions. Moreover, in the example they implemented which is based on a steady advection-diffusion problem first suggested by \cite{Choi2005}, they noted that ``this choice of source function and kernel causes fortuitous algebraic simplifications." Therefore this method, whilst interesting as an iterative procedure, cannot be considered to have a general validity.

\subsection{Noble's matrix factorisation approach}
\label{sec:Noble}

To avoid the iterations, we tried to solve the two simultaneous Wiener-Hopf equations casted in matrix form according to the classic approach of \citet[pp.~153--157]{Noble1958}; see also \citet{Daniele1984} and \citet[Sec.~1.5.2]{Daniele2014}.
%The idea is to write the problem in the form
%\begin{equation}
%\label{eq:matrixWH1}
%\widehat{\mathbf{A}}\,\widehat{\mathbf{f}_+} + \widehat{\mathbf{B}}\,\widehat{\mathbf{f}_-} = \widehat{\mathbf{b}}
%\end{equation}
%and to transform it into
%\begin{equation}
%\label{eq:matrixWH2}
%\widehat{\mathbf{L}_+}\,\widehat{\mathbf{f}_+} + \widehat{\mathbf{L}_-}\,\widehat{\mathbf{f}_-} = \widehat{\mathbf{c}},
%\end{equation}
%after which it can be solved decomposing the right-hand side $\widehat{\mathbf{c}}$. The step from Eq.~(\ref{eq:matrixWH1}) to Eq.~(\ref{eq:matrixWH1}) can be achieved multiplying Eq.~(\ref{eq:matrixWH1}) from the left by $\widehat{\mathbf{A}}^{-1}$ or $\widehat{\mathbf{B}}^{-1}$, factorising $\widehat{\mathbf{B}}^{-1}\widehat{\mathbf{A}} =: \widehat{\mathbf{L}} = \widehat{\mathbf{L}}_-\widehat{\mathbf{L}}_+$, and multiplying from the left by $\widehat{\mathbf{L}}_-^{-1}$.
We write Eq.~(\ref{eq:F_FT}) multiplied once by $e^{-ia\xi}$ and once by $e^{-ib\xi}$ as
\begin{equation}
\begin{pmatrix} \widehat{l} & e^{i(d-a)\xi} \\ 0 & e^{i(d-b)\xi} \end{pmatrix}
\begin{pmatrix} \widehat{f_0}e^{-ia\xi} \\ \widehat{f_+}e^{-id\xi} \end{pmatrix} +
\begin{pmatrix} 0 & e^{i(c-a)\xi} \\ \widehat{l} & e^{i(c-b)\xi} \end{pmatrix}
\begin{pmatrix} \widehat{f_0}e^{-ib\xi} \\ \widehat{f_-}e^{-ic\xi} \end{pmatrix} =
\begin{pmatrix} 0 & 1 \\ \widehat{l} & e^{i(a-b)\xi} \end{pmatrix}
\begin{pmatrix} 0 \\ \widehat{g_0}e^{-ia\xi} \end{pmatrix}
\end{equation}
where $a \le c$ and $d \le b$. Convenient choices of the parameters $c$ and $d$ are $c = a,\ d = b$; $c = b,\ d = a$; $c = d = a$; $c = d = b$. We choose $c = a,\ d = b$ and write for short
\begin{equation}
\label{eq:F_FT_N1}
\widehat{\mathbf{L}_1}\,\widehat{\mathbf{f}_+} + \widehat{\mathbf{L}_2}\,\widehat{\mathbf{f}_-} = \widehat{\mathbf{L}_2}\,\widehat{\mathbf{g_+}}.
\end{equation}
Multiplying from the left with $\widehat{\mathbf{L}_2}^{-1}$ yields a matrix version of Eq.~(\ref{eq:WH1_FT}),
\begin{equation}
\label{eq:F_FT_N2}
\widehat{\mathbf{L}}\,\widehat{\mathbf{f}_+} + \widehat{\mathbf{f}_-} = \widehat{\mathbf{g}_+},
\end{equation}
where
\begin{equation}
\widehat{\mathbf{L}} = \widehat{\mathbf{L}_2}^{-1}\,\widehat{\mathbf{L}_1} = \widehat{l}^{-1} \begin{pmatrix} -e^{i(a-b)\xi} & 1 \\ \widehat{l} & 0 \end{pmatrix} \begin{pmatrix} \widehat{l} & e^{i(b-a)\xi} \\ 0 & 1 \end{pmatrix} = \begin{pmatrix} -e^{i(a-b)\xi} & 0 \\ \widehat{l} & e^{i(b-a)\xi} \end{pmatrix}
\end{equation}
is a triangular matrix. Swapping the elements of $\widehat{\mathbf{f}_+}$ and $\widehat{\mathbf{f}_-}$ permutes the elements of $\widehat{\mathbf{L}}$. If we knew how to factorise $\widehat{\mathbf{L}} = \widehat{\mathbf{L}_-}\,\widehat{\mathbf{L}_+}$, multiplying Eq.~(\ref{eq:F_FT_N2}) from the left with $\widehat{\mathbf{L}_-}^{-1}$ would lead finally to
\begin{equation}
\label{eq:F_FT_N3}
\widehat{\mathbf{L}_+}\,\widehat{\mathbf{f}_+} + \widehat{\mathbf{L}_-}^{-1}\,\widehat{\mathbf{f}_-} = \widehat{\mathbf{L}_-}^{-1}\,\widehat{\mathbf{g}_+},
\end{equation}
which is a matrix version of Eq.~(\ref{eq:WH1_FT_2}) and can be solved in a similar fashion decomposing its right-hand side. The same result is obtained multiplying Eq.~(\ref{eq:F_FT_N1}) from the left with $\widehat{\mathbf{L}_1}^{-1}$ or Eq.~(\ref{eq:F_FT_N2}) from the left with $\widehat{\mathbf{L}}^{-1}$, yielding
\begin{equation}
\label{eq:F_FT_N4}
\widehat{\mathbf{f}_+} + \widehat{\mathbf{L}}^{-1}\,\widehat{\mathbf{f}_-} = \widehat{\mathbf{L}}^{-1}\,\widehat{\mathbf{g}_+},
\end{equation}
where
\begin{equation}
\widehat{\mathbf{L}}^{-1} = \widehat{\mathbf{L}_1}^{-1}\,\widehat{\mathbf{L}_2} = \widehat{l}^{-1} \begin{pmatrix} 1 & -e^{i(b-a)\xi} \\ 0 & \widehat{l} \end{pmatrix} \begin{pmatrix} 0 & 1 \\ \widehat{l} & e^{i(a-b)\xi} \end{pmatrix} = \begin{pmatrix} -e^{i(b-a)\xi} & 0 \\ \widehat{l} & e^{i(a-b)\xi} \end{pmatrix}
\end{equation}
If we knew how to factorise $\widehat{\mathbf{L}}^{-1} = \widehat{\mathbf{L}_+}^{-1}\,\widehat{\mathbf{L}_-}^{-1}$, multiplying Eq.~(\ref{eq:F_FT_N4}) from the left with $\widehat{\mathbf{L}_+}$ would lead again to Eq.~(\ref{eq:F_FT_N3}).

Unfortunately, though a formula to factorise triangular $2 \times 2$ matrices has been published, see \citet[Eq.~(21)]{Jones1984} and \citet[Eq.~(6)]{Jones1991}, an oscillatory element like $e^{i(b-a)\xi}$ does not fulfil the conditions to apply it: it is assumed that $+$ or $-$ factors remain $+$ or $-$ when inverted, but the inverse of the $+$ function $e^{i(b-a)\xi}$ is the $-$ function $e^{i(a-b)\xi}$; see also \citet[Sec.~4.3, Example 2]{Daniele2014}.

\subsection{Voronin's matrix factorisation approach}
\label{sec:Voronin}

More recently, a different matrix form of the two simultaneous Wiener-Hopf equations has been proposed by \cite{Voronin2004}. We present it with slight modifications. Start from Eq.~(\ref{eq:F_FT}) and decompose the kernel, $\widehat{k_0} = \widehat{k_-} + \widehat{k_+}$ (for simplicity we drop the 0 subscript from the right-hand side), obtaining
\begin{equation}
(\lambda-\widehat{k_-}-\widehat{k_+})\widehat{f_0} + \widehat{f_-} + \widehat{f_+} = \widehat{g_0}.
\end{equation}
Multiply by $e^{-ia\xi}$, take the $+$ part, thus eliminating $\widehat{f_-}$, which is a $-$ function with respect to $a$, and multiply by $e^{ia\xi}$, yielding
\begin{equation}
(\lambda-\widehat{k_+})\widehat{f_0} - e^{ia\xi}(\widehat{k_-}\widehat{f_0}e^{-ia\xi})_+ + \widehat{f_+} = \widehat{g_0}.
\end{equation}
Multiply by $e^{-ib\xi}$, take the $-$ part, thus eliminating $\widehat{f_+}$, which is a $+$ function with respect to $b$, and multiply by $e^{ib\xi}$, yielding
\begin{equation}
\lambda\widehat{f_0} - e^{ia\xi}(\widehat{k_-}\widehat{f_0}e^{-ia\xi})_+ - e^{ib\xi}(\widehat{k_+}\widehat{f_0}e^{-ib\xi})_- = \widehat{g_0}.
\end{equation}
Define $\widehat{\varphi_1} = \big(\widehat{k_-}+\frac{1-\lambda}{2}\big)\widehat{f_0}$ and decompose it with respect to $a$, $\widehat{\varphi_1} = \widehat{\varphi_{1+}} + \widehat{\varphi_{1-}}$; define $\widehat{\varphi_2} = \big(\widehat{k_+}+\frac{1-\lambda}{2}\big)\widehat{f_0}$ and decompose it with respect to $b$, $\widehat{\varphi_2} = \widehat{\varphi_{2+}} + \widehat{\varphi_{2-}}$; this gives
\begin{equation}
\label{eq:F_FT_V1}
\widehat{f_0} - \widehat{\varphi_{1+}} - \widehat{\varphi_{2-}} = \widehat{g_0}.
\end{equation}
The two coupled Wiener-Hopf equations are now obtained multiplying once by $\widehat{k_-}e^{-ia\xi}$ and once by $\widehat{k_+}e^{-ib\xi}$,
\begin{equation}
\begin{pmatrix} 1-\widehat{k_-} & 0 \\ -\widehat{k_+}e^{i(a-b)\xi} & 1 \end{pmatrix}
\begin{pmatrix} \widehat{\varphi_{1+}}e^{-ia\xi} \\ \widehat{\varphi_{2+}}e^{-ib\xi} \end{pmatrix} +
\begin{pmatrix} 1 & -\widehat{k_-}e^{i(b-a)\xi} \\ 0 & 1-\widehat{k_+} \end{pmatrix}
\begin{pmatrix} \widehat{\varphi_{1-}}e^{-ia\xi} \\ \widehat{\varphi_{2-}}e^{-ib\xi} \end{pmatrix} =
\begin{pmatrix} \widehat{k_-}\widehat{g_0}e^{-ia\xi} \\ \widehat{k_+}\widehat{g_0}e^{-ib\xi} \end{pmatrix},
\end{equation}
for short
\begin{equation}
\label{eq:F_FT_V2}
\widehat{\mathbf{M}_{\mathrm{r}-}}^{-1}\widehat{\boldsymbol{\varphi}_+} + \widehat{\mathbf{M}_{\mathrm{r}+}}\,\widehat{\boldsymbol{\varphi}_-} = \widehat{\mathbf{g}}.
\end{equation}
Here one can see that the parameter $\lambda$ has been inserted in the definition of $\widehat{\boldsymbol{\varphi}_1}$ and $\widehat{\boldsymbol{\varphi}_2}$ to avoid that it appears in place of the numbers 1 in the diagonal elements of $\widehat{\mathbf{M}_{\mathrm{r}-}}^{-1}$ and $\widehat{\mathbf{M}_{\mathrm{r}+}}$, which would make these matrices singular for $\lambda = 0$. Multiplying from the left by $\widehat{\mathbf{M}_{\mathrm{r}-}}$ yields
\begin{equation}
\label{eq:F_FT_V3}
\widehat{\boldsymbol{\varphi}_+} + \widehat{\mathbf{M}}\,\widehat{\boldsymbol{\varphi}_-} = \widehat{\mathbf{M}_{\mathrm{r}-}}\,\widehat{\mathbf{g}},
\end{equation}
where
\begin{eqnarray}
\label{eq:M-+}
\widehat{\mathbf{M}} &=& \widehat{\mathbf{M}_{\mathrm{r}-}}\widehat{\mathbf{M}_{\mathrm{r}+}}
= \frac{1}{1-\widehat{k_-}}
\begin{pmatrix} 1 & 0 \\ \widehat{k_+}e^{i(a-b)\xi} & 1-\widehat{k_-} \end{pmatrix}
\begin{pmatrix} 1 & -\widehat{k_-}e^{i(b-a)\xi} \\ 0 & 1-\widehat{k_+} \end{pmatrix} \nonumber \\
&=& \frac{1}{1-\widehat{k_-}}
\begin{pmatrix} 1 & -\widehat{k_-}e^{i(b-a)\xi} \\ \widehat{k_+}e^{i(a-b)\xi} & 1-\widehat{k} \end{pmatrix}.
\end{eqnarray}
An equivalent result is obtained multiplying Eq.~(\ref{eq:F_FT_V2}) from the left by $\widehat{\mathbf{M}_{\mathrm{r}+}}^{-1}$. If we knew how to factorise $\widehat{\mathbf{M}} = \widehat{\mathbf{M}_{\mathrm{l}+}}\widehat{\mathbf{M}_{\mathrm{l}-}}$, multiplying Eq.~(\ref{eq:F_FT_V3}) from the left by $\widehat{\mathbf{M}_{\mathrm{l}+}}^{-1}$ would lead finally to
\begin{equation}
\widehat{\mathbf{M}_{\mathrm{l}+}}^{-1}\widehat{\boldsymbol{\varphi}_+} + \widehat{\mathbf{M}_{\mathrm{l}-}}\,\widehat{\boldsymbol{\varphi}_-} = \widehat{\mathbf{M}_{\mathrm{l}+}}^{-1}\widehat{\mathbf{M}_{\mathrm{r}-}}\,\widehat{\mathbf{g}},
\end{equation}
which, like Eq.~(\ref{eq:F_FT_N3}), is a matrix version of Eq.~(\ref{eq:WH1_FT_2}) and can be solved in a similar fashion decomposing its right-hand side.

Unfortunately we are stuck again: though formulas to convert left ($+-$) factorisations of $2 \times 2$ matrices into right ($-+$) ones and vice versa have been published by \citet[Eqs.~(8)--(11)]{Jones1991}, the information contained in Eq.~(\ref{eq:M-+}) is not sufficient to apply them, so we cannot obtain $\widehat{\mathbf{M}_{\mathrm{l}+}}$ and $\widehat{\mathbf{M}_{\mathrm{l}-}}$ from our knowledge of $\widehat{\mathbf{M}_{\mathrm{r}-}}$ and $\widehat{\mathbf{M}_{\mathrm{r}+}}$.

\subsection{Iterative solution based on Voronin's approach}
\label{sec:iterative_Voronin}

An iterative solution is possible also with Voronin's approach. Write Eq.~(\ref{eq:F_FT_V1}) multiplied once by $\widehat{k_-}/(1-\widehat{k_-})$ and once by $\widehat{k_+}/(1-\widehat{k_+})$,
\begin{eqnarray}
\widehat{\varphi_{1+}} - \frac{1}{1-\widehat{k_-}}\widehat{\varphi_{1-}} - \frac{\widehat{k_-}}{1-\widehat{k_-}}\widehat{\varphi_{2-}} &=& \frac{\widehat{k_-}}{1-\widehat{k_-}}\widehat{g_0} \label{eq:F_FT_iV1} \\
\widehat{\varphi_{2-}} - \frac{1}{1-\widehat{k_+}}\widehat{\varphi_{2+}} - \frac{\widehat{k_+}}{1-\widehat{k_+}}\widehat{\varphi_{1+}} &=& \frac{\widehat{k_+}}{1-\widehat{k_+}}\widehat{g_0}, \label{eq:F_FT_iV2}
\end{eqnarray}
in Eq.~(\ref{eq:F_FT_iV1}) set
\begin{equation}
\widehat{c_1} = \frac{\widehat{k_-}}{1-\widehat{k_-}}(\widehat{g_0}+\widehat{\varphi_{2-}}),
\end{equation}
start from the guess $\widehat{\varphi_{2-}} = 0$, decompose $\widehat{c_1} = \widehat{c_{1+}} + \widehat{c_{1-}}$ with respect to $a$, and compute the approximations
\begin{gather}
\widehat{\varphi_{1+}} = \widehat{c_{1+}}\\
\widehat{\varphi_{1-}} = (1-\widehat{k_-})\widehat{c_{1-}}
\end{gather}
as $+$ and $-$ functions with respect to $a$, as well as
\begin{equation}
\widehat{f_0} = \frac{1}{\widehat{k_-}+\frac{1-\lambda}{2}}(\widehat{\varphi_{1+}} + \widehat{\varphi_{1-}}) = \frac{1}{\widehat{k_-}+\frac{1-\lambda}{2}}\big(\widehat{c_{1+}} + (1-\widehat{k_-})\widehat{c_{1-}}\big);
\end{equation}
then turn to Eq.~(\ref{eq:F_FT_iV2}), set
\begin{equation}
\widehat{c_2} = \frac{\widehat{k_+}}{1-\widehat{k_+}}(\widehat{g_0}-\widehat{\varphi_{1+}}),
\end{equation}
decompose $\widehat{c_2} = \widehat{c_{2+}} + \widehat{c_{2-}}$ with respect to $b$, and compute the new approximations
\begin{gather}
\widehat{\varphi_{2-}} = \widehat{c_{2-}}\\
\widehat{\varphi_{2+}} = (1-\widehat{k_+})\widehat{c_{2+}}
\end{gather}
as $+$ and $-$ functions with respect to $b$, as well as
\begin{equation}
\widehat{f_0} = \frac{1}{\widehat{k_+}+\frac{1-\lambda}{2}}(\widehat{\varphi_{2-}} + \widehat{\varphi_{2+}}) = \frac{1}{\widehat{k_+}+\frac{1-\lambda}{2}}\big(\widehat{c_{2-}} + (1-\widehat{k_+})\widehat{c_{2+}}\big);
\end{equation}
and so on until the difference between the values of $\widehat{f_0}$ at an iteration and the previous falls below a threshold. An equivalent result is obtained starting from the guess $\widehat{\varphi_{1+}} = 0$ in Eq.~(\ref{eq:F_FT_iV2}) and the computation of $\widehat{c_2}$. Similarly to Sec.~\ref{sec:iterative}, the iterations are performed looking for a fixed point on the variables $\widehat{\varphi_1}$ and $\widehat{\varphi_2}$, while $\widehat{f_0}$ is a side product output at each step, but not used to compute the next step. For details, see the function whf\_gmf\_v.m in the supplementary material. % Sec.~\ref{sec:whf_gmf_v.m}.

%\textcolor{red}{\textit{The idea was that Voronin's auxiliary variables $\widehat{\varphi_1}$ and $\widehat{\varphi_2}$ might be numerically more convenient than $\widehat{f_-}$ and $\widehat{f_+}$. Unfortunately this works only in the simple numerical tests below, but not when used to price double-barrier options. Moreover oddly I was not able to recover a non-iterative solution of the CWHE for $\lim_{a \to -\infty}$ or $\lim_{b \to -\infty}$, or actually $a = -t$ or $b = t$.}}

Generally, work on factorisation has concentrated on finding methods for particular matrix classes, often corresponding to specific applications. For example, \cite{Kisil2015} has developped and analysed an approximate factorisation approach for Daniele-Khrapkov matrics but, although several matrix classes can be reduced to this case, this is not a general solution. A digest of constructive matrix factorisation methods was published in this journal by \cite{Rogosin2016}, who described several different matrix classes and the most suitable methods for their factorisation. This work, whilst important, underlines that there is no general method which can be used for all matrices.

%we define more accurately the iterative refinement and evaluate its performance

\section{Test cases}\label{sec:8_Test}
%\textcolor{red}{Explain in a couple of lines why we present 3 test cases for the Fredholm equation and none for the Wiener-Hopf equation: Fredholm is more challenging and contains Wiener-Hopf as a special case when $a \to -\infty$ or $b \to +\infty$. Actually in a numerical scheme both extremes are always finite, but if one is large enough in absolute value, the algorithm ceases to be iterative. A few similar sentences must be added also at the beginning of the Conclusion section, which mentions only the Fredholm equation although the paper is also about the Wiener-Hopf equation, i.e.\ on a semi-infinite domain. Eliminate Eq. (5.1) and instead refer to Eq.~(\ref{eq:WHF}) which is identical.}

As we present a general solution to the Fredholm equation, rather than one limited to a particular application, we provide several test cases to solve Eq.~(\ref{eq:WHF})
%\begin{align}
%\lambda f(x)-\int^b_a k(x-x')f(x')dx'=g(x), \quad x\in(a,b)\label{eq:8_an_whf}
%\end{align}
for $f(x)$. %As we do not have a numerical method with very high accuracy that can be used as a comparison (as we did with the method by \cite{Feng2008} for barrier pricing applications) we require test cases with analytic solutions.
Although the methods developped herein can be applied to both Fredholm and Wiener-Hopf equations, for the numerical tests we have chosen to concentrate on solving examples of the Fredholm equation as it is the more challenging case and encompasses Wiener-Hopf as a special example when $a \to -\infty$ or $b \to +\infty$.

Solutions to Eq.~(\ref{eq:WHF}) with simple closed form expressions for $f(x)$, $g(x)$ and $k(x)$ are not readily available. However, if we limit the requirement for simplicity to $f(x)$ and $k(x)$, then closed form solutions for $g(x)$ in Eq.~(\ref{eq:WHF}) can be calculated. %The solution methods can then be compared by applying them to $g(x)$ and $k(x)$, solving numerically for $f(x)$ and calculating the error compared to the analytic expression.
These $g(x)$ and $k(x)$ can be used as inputs to our numerical method, whose accuracy can be measured by comparing the result with $f(x)$.
We selected $f(x)$ and $k(x)$ to give closed form expressions for $g(x)$ and also to have Fourier transforms which are easily calculable. We derived three solutions, with both $f(x)$ and $k(x)$ 1. Gaussian, 2. Cauchy, 3. Laplace (bilateral exponential). As discussed in Section \ref{sec:Gibbs_phenom} the decay of the functions as $x,\xi\rightarrow\infty$ can influence the error performance of Fourier-based methods. The functions were therefore selected to be exponentially decaying in both the state space and Fourier space (Gaussian), to be polynomially decaying in the state space and exponentially decaying in the Fourier space (Cauchy) and to be exponentially decaying in the state space and polynomially decaying in the Fourier space (Laplace). The derivation of $g(x)$ is described in the following sections. 

\subsection{Gaussian}\label{sec:8_Test_Gauss}
We set $f(x)=k(x)=\frac{1}{\sqrt{\pi}}e^{-x^2}$. The expression for $g(x)$ is then
\begin{align}
g(x)&=\frac{1}{\sqrt{\pi}}e^{-x^2}-\frac{1}{\pi}\int^b_ae^{-(x-y)^2}e^{-y^2}dy\nonumber\\
&=\frac{1}{\sqrt{\pi}}e^{-x^2}-\frac{1}{\pi}e^{-\frac{x^2}{2}}\int^b_ae^{-2(y-\frac{x}{2})^2}dy\nonumber\\
&=\frac{1}{\sqrt{\pi}}e^{-x^2}-\frac{1}{\sqrt{2\pi}}e^{-\frac{x^2}{2}}\big[\Phi(2b-x)-\Phi(2a-x)\big],\quad x\in[a,b],
\end{align}
where $\Phi(\cdot)$ is the standard normal cumulative distribution function.

\subsection{Cauchy}\label{sec:8_Test_Cauch}
We set $f(x)=k(x)=\frac{1}{\pi(x^2+1)}$. The first step in the calculation of $g(x)$ is to solve the integral
\begin{align}
g_{\mathrm{int}}(x) & =\frac{1}{\pi^2}\int^b_a\frac{1}{y^2+1}\frac{1}{(x-y)^2+1}dy.\label{eq:fred}
\end{align}
Using partial fractions,
\begin{align}
g_{\mathrm{int}}(x) =&\frac{1}{\pi^2}\int^b_a\frac{1}{y^2+1}\frac{1}{(x-y)^2+1}dy\nonumber\\
%=&\frac{1}{\pi^2x(x^2+4)}\int^b_a\frac{2y}{y^2+1}+\frac{x}{y^2+1}-\frac{2y}{(y-x)^2+1}+\frac{3x}{(y-x)^2+1}dy\nonumber\\
=&\frac{1}{\pi^2x(x^2+4)}\int^b_a\left(\frac{2y}{y^2+1}+\frac{x}{y^2+1}-\frac{2(y-x)}{(y-x)^2+1}+\frac{x}{(y-x)^2+1}\right)dy\nonumber\\
=&\frac{1}{\pi^2x(x^2+4)}\left[\log(y^2+1)+x\arctan(y)-\log[(x-y)^2+1]+x\arctan(y-x)\right]^b_a\nonumber\\
=&\frac{1}{\pi^2x(x^2+4)}\left\{\log\frac{(b^2+1)((a-x)^2+1)}{(a^2+1)((b-x)^2+1)}+\right.\nonumber\\
&\left.+x\left[\arctan(b)-\arctan(a)+\arctan(b-x)-\arctan(a-x)\right]\right\}.
\end{align}
This gives $g(x)$ in closed form:
\begin{align}
g(x)=&\frac{1}{\pi(x^2+1)}-\frac{1}{\pi^2x(x^2+4)}\left\{\log\frac{(b^2+1)((a-x)^2+1)}{(a^2+1)((b-x)^2+1)}+\right.\nonumber\\
&\left.+x\left[\arctan(b)-\arctan(a)+\arctan(b-x)-\arctan(a-x)\right]\right\}, \quad x\in[a,b].\label{eq:cauchg}
\end{align}

\subsection{Laplace or bilateral exponential}\label{sec:8_Test_Exp}
We set $f(x)=k(x)=\frac{1}{2}e^{-|x|}$. In order to make the calculation of $g(x)$ simpler, the values of $a$ and $b$ are restricted so that $0<a<b<\infty$. Then the formula for $g(x)$ in closed form is
\begin{align}
g(x)&=\frac{1}{2}e^{-x}-\frac{1}{4}\int^b_ae^{-|x-y|}e^{-y}dy\nonumber\\
&=\frac{1}{2}e^{-x}-\frac{1}{4}\left[\int^b_xe^{(x-y)}e^{-y} dy +\int^x_a e^{-(x-y)}e^{-y} dy\right]\nonumber\\
&=\frac{1}{2}e^{-x}-\frac{1}{4}e^x\left(\int^b_xe^{-2y} dy +e^{-x}\int^x_a dy\right)\nonumber\\
&=\frac{1}{2}e^{-x}+\frac{1}{8}e^x\left[e^{-2y}\right]^b_x-\frac{1}{4}e^{-x}\big[y\big]^x_a\nonumber\\
&=\frac{1}{2}e^{-x}+\frac{1}{8}e^x\left(e^{-2b}-e^{-2x}\right)-\frac{1}{4}e^{-x}(x-a)\nonumber\\
&=e^{-x}\left[\frac{3}{8}+\frac{1}{8}e^{-2(b-x)}+\frac{1}{4}(a-x)\right], \quad x\in[a,b].
\end{align}

\section{Results}\label{sec:8_Res}
We used the following methods to recover $f(x)$ and produce the detailed results shown in this section:
\begin{enumerate}
\item 4th order quadrature with preconditioner \citep{Fusai2012,Press2007}; see the \textsc{Matlab} functions quadrature.m and weights.m in the supplementary material.
\item Wiener-Hopf method using the sinc-based fast Hilbert transform with no zero padding. In order to counteract the oscillations on the recovered function, we used an exponential filter of order 8 on the final stage of the fixed point algorithm. The maximum number of iterations of the fixed point algorithm is set to 5. In fact, in most cases the final error level is achieved within 3 iterations. We discuss the use of the sinc-based fast Hilbert transform and spectral filtering in Section \ref{sec:8_sincfht} below.
\item Wiener-Hopf method using the symmetric sign function in the fast Hilbert transform, i.e.\ with zeros placed at both $\xi=0$ and $\xi=\xi_{\min}=-\frac{N}{2}\Delta\xi$, similar to the method introduced by \cite{Rino1970} and \cite{Henery1974} and tested by \cite{Fusai2016}.
\item Wiener-Hopf method with Voronin's variant using the symmetrical sign function for the Hilbert transform.
\end{enumerate}
%It is common in the literature on numerical methods to refer to the number of grid points as $n$ or $N$, however to maintain consistency with the option pricing methods described in Chapters \ref{chapterlabel4}--\ref{chapterlabel7}, where $N$ is used for the number of monitoring dates, we use $M$.
\subsection{Results for the Gaussian test case}
We first examine the performance of the different numerical methods with the Gaussian test case, with particular emphasis on the method used to implement the Hilbert transform.
\subsubsection{Sinc-based fast Hilbert transform and spectral filtering}\label{sec:8_sincfht}
In the financial pricing applications described by \cite{Feng2008} and \cite{Fusai2016} the sinc-based fast Hilbert transform has shown excellent error convergence, especially when combined with a spectral filter as in \cite{Phelan2019}. However, when we consider its use for this application we must take account of several ways in which the requirements differ from its general use for finding solutions to Wiener-Hopf or Fredholm equations.

Firstly the pricing methods that were implemented using the Wiener-Hopf method in \cite{Fusai2016}, as devised by \cite{Green2010}, use the analytic continuation of $x$, i.e.\ they give results for values of $x$ both inside and outside the barriers (the integration limits of Fredholm equation). This means that for these applications there is no requirement to truncate the functions to the integration limits of $a$ and $b$ in the state space, unlike the problems presented as examples in this paper. The requirement to truncate the function means that there is a jump discontinuity introduced in the state space, meaning that, as described by \cite{Boyd2001}, the function in the Fourier space decays as a first order polynomial due to the Gibbs phenomenon. As explained in \cite{Stenger1993} this polynomial decay means that the sinc-based fast Hilbert transform no longer obtains an error which is exponentially convergent with grid size but rather converges polynomially. This is in contrast with the aforementioned finance-based papers from \cite{Green2010} and \cite{Fusai2016} where the Fourier domain functions subject to the sinc-based fast Hilbert transform are exponentially decaying (or polynomially so in the case of the VG process) and so excellent error performance is achieved, especially in conjunction with a spectral filter to solve the issue with the fixed point algorithm.

In contrast, here we solve the Fredholm and Wiener-Hopf equations as they were originally formulated, i.e.\ the function is only defined for the range of the integration $[a,b]$ and therefore the functions $g(x)$ and $k(x)$ must be truncated to the ranges $[a,b]$ and $[a-b,b-a]$ respectively. This truncation will introduce a jump in the functions which means that their Fourier transforms now have first order polynomial decay. Therefore the truncation error from the Hilbert transform will have a first order polynomial convergence unless we can exploit some symmetry between the Fourier domain functions for positive and negative $\xi$ as in \cite{Phelan2018}, in which case we may achieve second order polynomial convergence.

Moreover, there is a second important distinction to be made between the general solution presented here and the work in the above literature. In the finance literature the solutions to the Fredholm equation are used to calculate the expectation of a further function, in this case the payoff function. Therefore the exact errors in the function for individual values of $x$ are not particularly important. Rather, the finance literature is concerned with the average error, weighted according to the shape of the payoff function. This also has particular importance when we are considering the use of the sinc-based fast Hilbert transform described in Section \ref{sec:HilbertSinc} which was instrumental in achieving exponential error convergence with the number of FFT grid points $N$ in \cite{Feng2008} and \cite{Fusai2016}.

%A second difference when we are using this numerical method for pricing is that, after the Fredholm equation is solved to provide the Sptizer identities, we then multiply the result by the payoff function in the Fourier domain which has a smoothing effect on the output. For the general numerical solution we have no such intrinsic smoothing, although we discuss the effects of smoothing below.
%
%The third major difference is that, for the pricing techniques, we use the Plancherel relation to extract the final price and thus are only interested in the solution at $x=0$, i.e.\ the centre of our state-space range. In the general case we are, of course, potentially interested in the solution at any point in the range $x\in(a,b)$.

\begin{figure}[h]
\vspace{-4mm}
\begin{center}
\includegraphics[width=\textwidth]{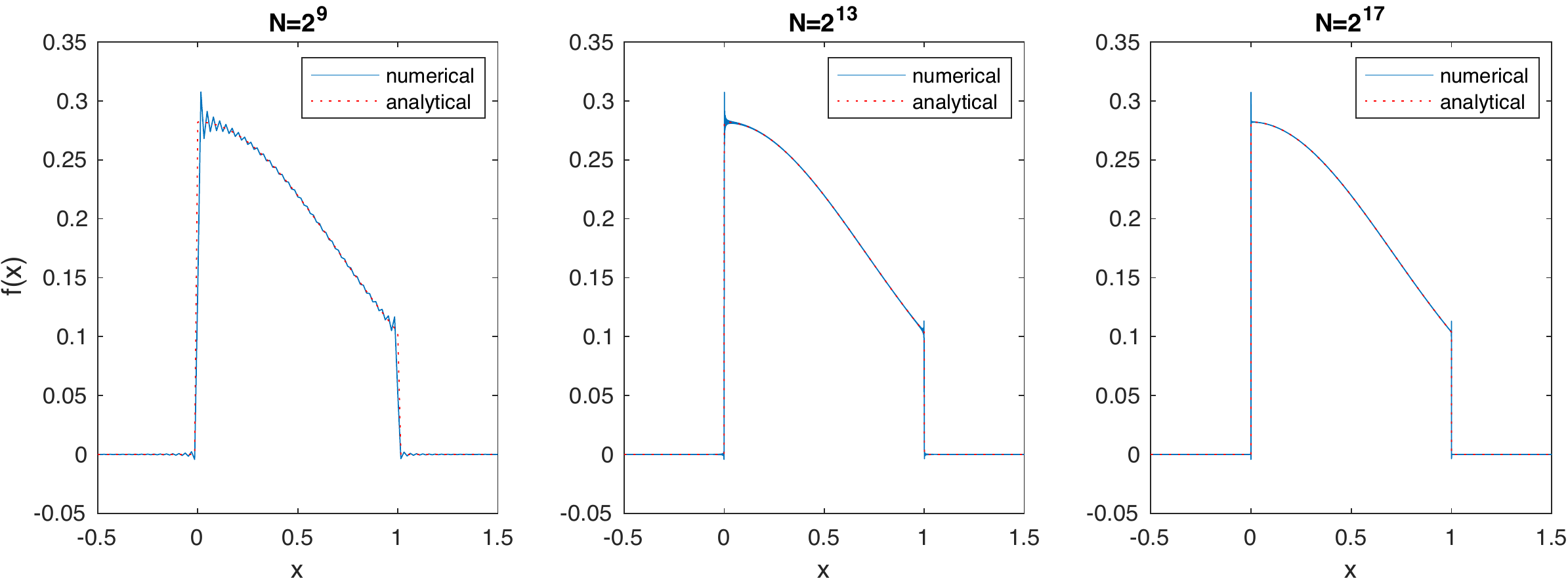}
\caption{Numerical and analytical $f(x)$ using the sinc-based fast Hilbert transform with no filter.}
\label{fig:gaussinc}
\end{center}
\vspace{-4mm}
\end{figure}

\begin{figure}[h]
\vspace{-4mm}
\begin{center}
\includegraphics[width=\textwidth]{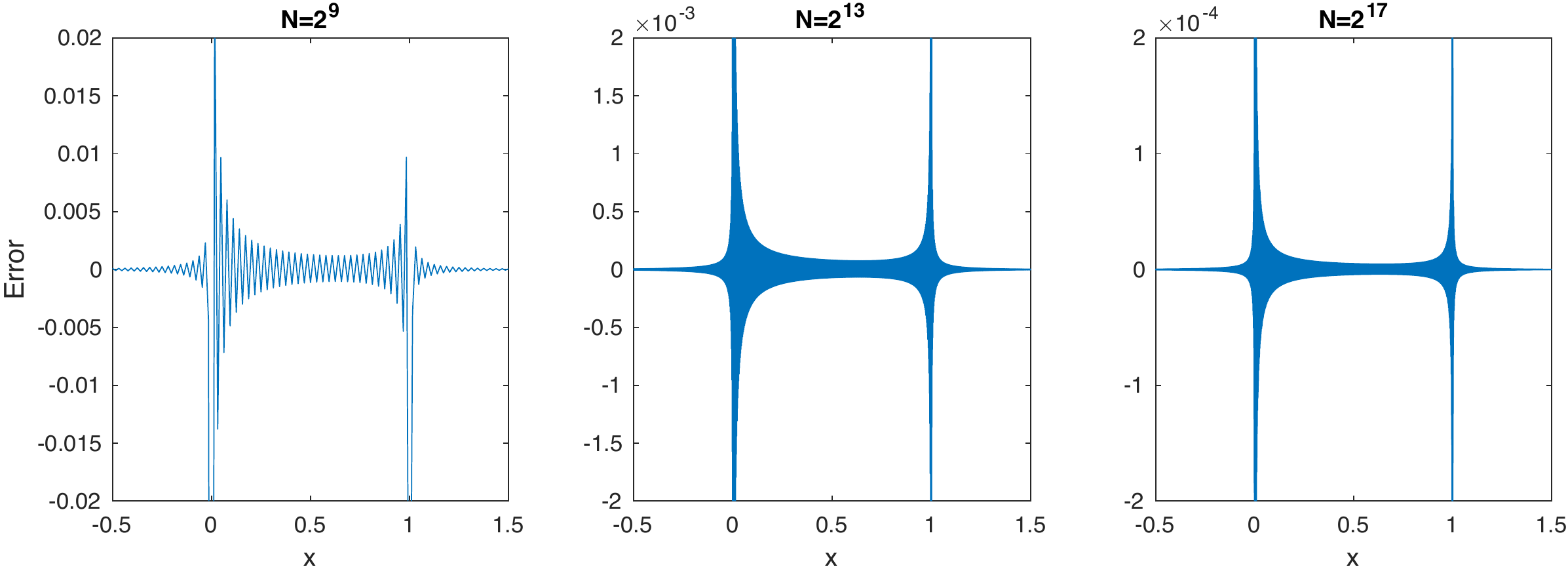}
\caption{Error in the numerical calculation of $f(x)$ using the sinc-based fast Hilbert transform with no filter.}
\label{fig:gaussincerr}
\end{center}
\vspace{-4mm}
\end{figure}

In Figures \ref{fig:gaussinc} and \ref{fig:gaussincerr}, we show results using the sinc-based fast Hilbert transform with no filter for the Gaussian test case described in Section \ref{sec:8_Test_Gauss}. It is immediately obvious that, even for high values of $N$, oscillations are visible in the numerical solution.

\begin{figure}[h]
\vspace{-4mm}
\begin{center}
\includegraphics[width=\textwidth]{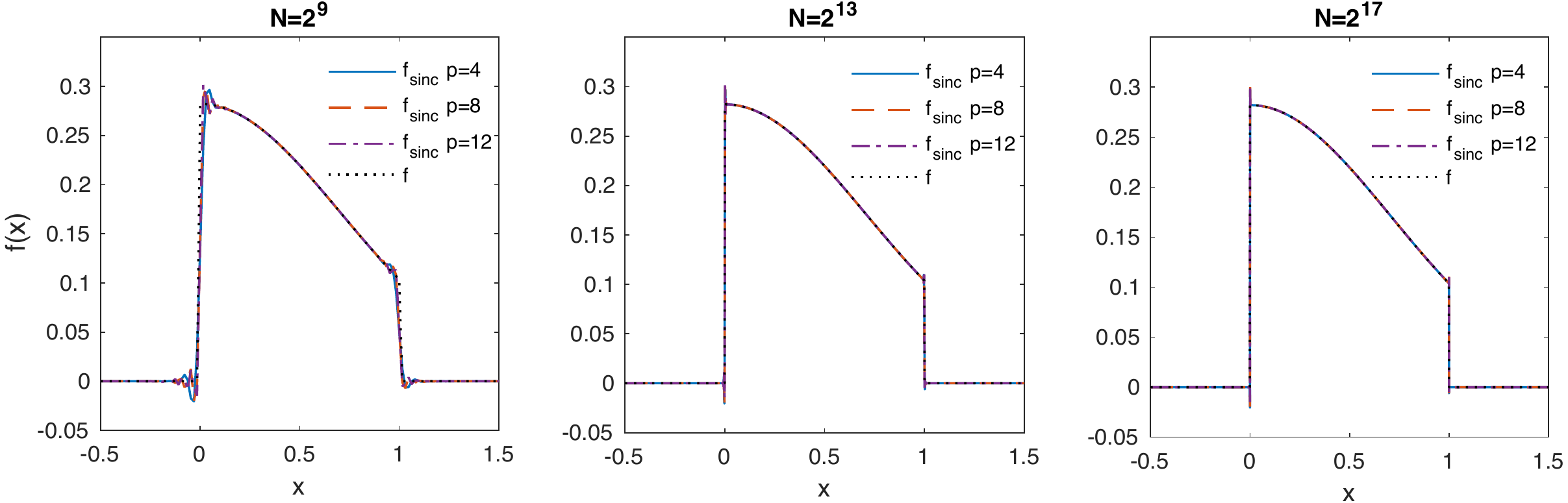}
\caption{Numerical and analytical $f(x)$ using the sinc-based fast Hilbert transform with an exponential filter. The parameter $p$ describes the order of the filter.}
\label{fig:gaussincfiltfull}
\end{center}
\vspace{-3mm}
\end{figure}

\begin{figure}[h]
\vspace{-4mm}
\begin{center}
\includegraphics[width=\textwidth]{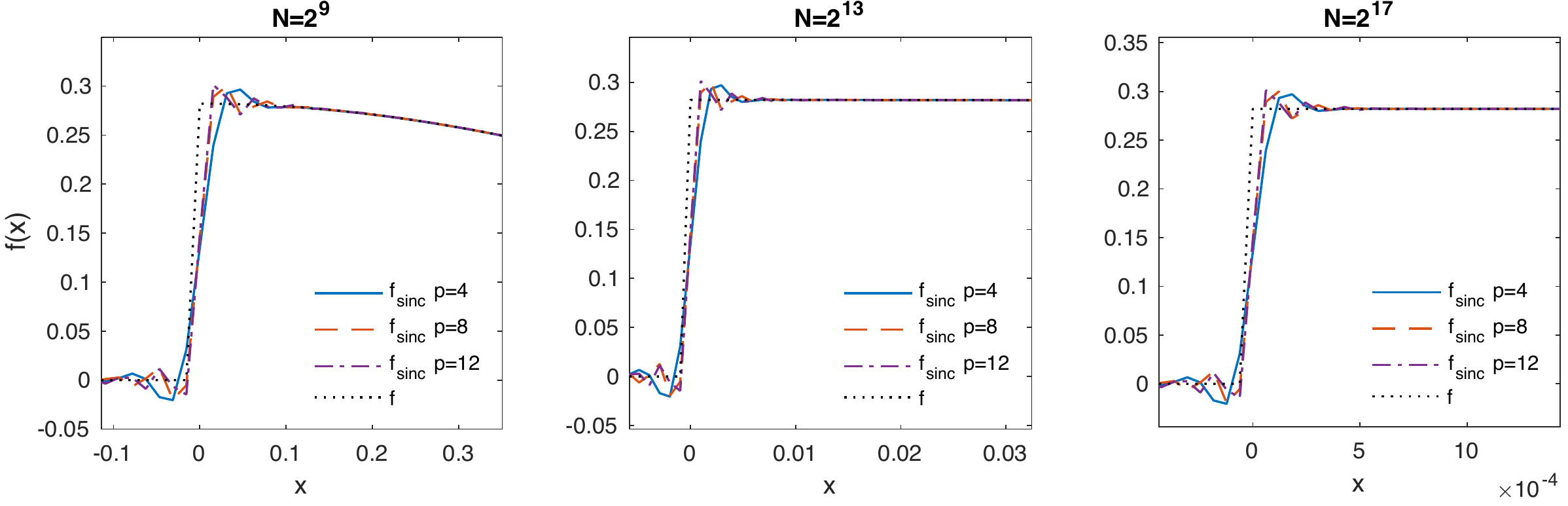}
\caption{Numerical and analytical $f(x)$ using the sinc-based fast Hilbert transform with an exponential filter, focussing on the discontinuity at $x=0$. The parameter $p$ describes the order of the filter.}
\label{fig:gaussincfiltover}
\end{center}
\vspace{-3mm}
\end{figure}

\begin{figure}
\vspace{-4mm}
\begin{center}
\includegraphics[width=\textwidth]{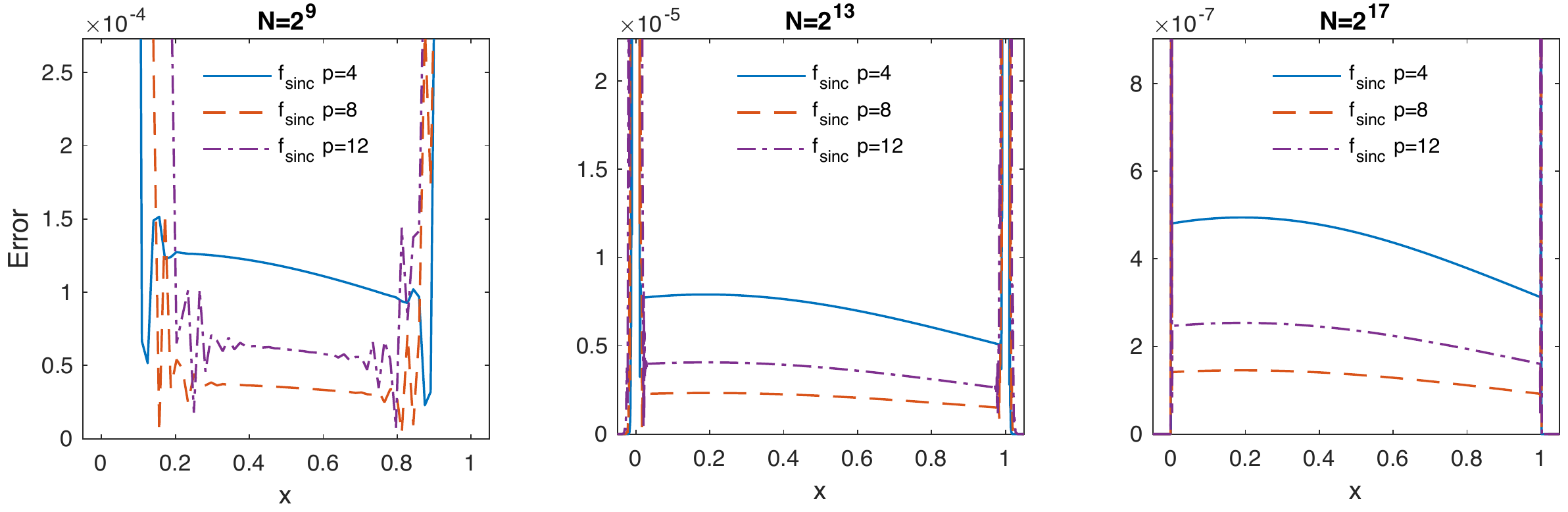}
\caption{Error between the numerical and analytical calculation of $f(x)$ using the sinc-based fast Hilbert transform with an exponential filter. The parameter $p$ describes the order of the filter. The scale has been chosen to display the error away from the discontinuities of $f(x)$.}
\label{fig:gaussincfilterror}
\end{center}
\vspace{-3mm}
\end{figure}

We can use a spectral filter to overcome the oscillations, however this can have a negative effect on the accuracy of the numerical method, especially close to the discontinuities in the state space; this is illustrated in Figures \ref{fig:gaussincfiltfull}--\ref{fig:gaussincfilterror}. Figure \ref{fig:gaussincfiltover} shows that the lower order filter gives a shallower slope at the discontinuity, but has a stronger effect on the oscillations. However, we can see from Figure \ref{fig:gaussincfiltover} that, regardless of the order of the filter, the overshoot at the discontinuity remains approximately the same. Figure \ref{fig:gaussincfilterror} shows that a spectral filter removes the oscillations away from the discontinuity and that the best results are achieved with a filter of order 8. Although the behaviour of the numerical method using the sinc-based fast Hilbert transform is not appropriate for a general solution to the Fredholm equation due to the high errors at function discontinuities, it remains the case that for applications where we are solely interested in a function value away from any jumps this may be an appropriate method to use.

\FloatBarrier
\subsubsection{Sign-based fast Hilbert transform method}
As an alternative to the sinc-based fast Hilbert transform, we examine the method of \cite{Rino1970} and \cite{Henery1974}, which was used also by \cite{Fusai2016}. It is based on the simple relationship between the Hilbert transform and the Fourier transform given in Eq.~(\ref{eq:HilbertFourier2}).%, i.e.
%\begin{align}
%\mathcal{H}[\widehat{f}(\xi)]=-i\mathcal{F}_{\xi\rightarrow x}[\mathrm{sgn}(x)\mathcal{F}^{-1}_{\xi\rightarrow x}\widehat{f}(\xi)],
%\end{align}
%where
%\begin{align}
%\mathrm{sgn}\,x=
%\begin{cases}
%	1\qquad &x>0 \\
%	0 \qquad &x=0 \\
%	-1 \quad &x<0 .
%\end{cases}\label{eq:8_sgn}
%	\end{align}
%Notice that Eq.~(\ref{eq:8_sgn}) is $0$ at $x=0$ and this is is reflected in our numerical implementation.

\begin{figure}[h]
\begin{center}
\includegraphics[width=\textwidth]{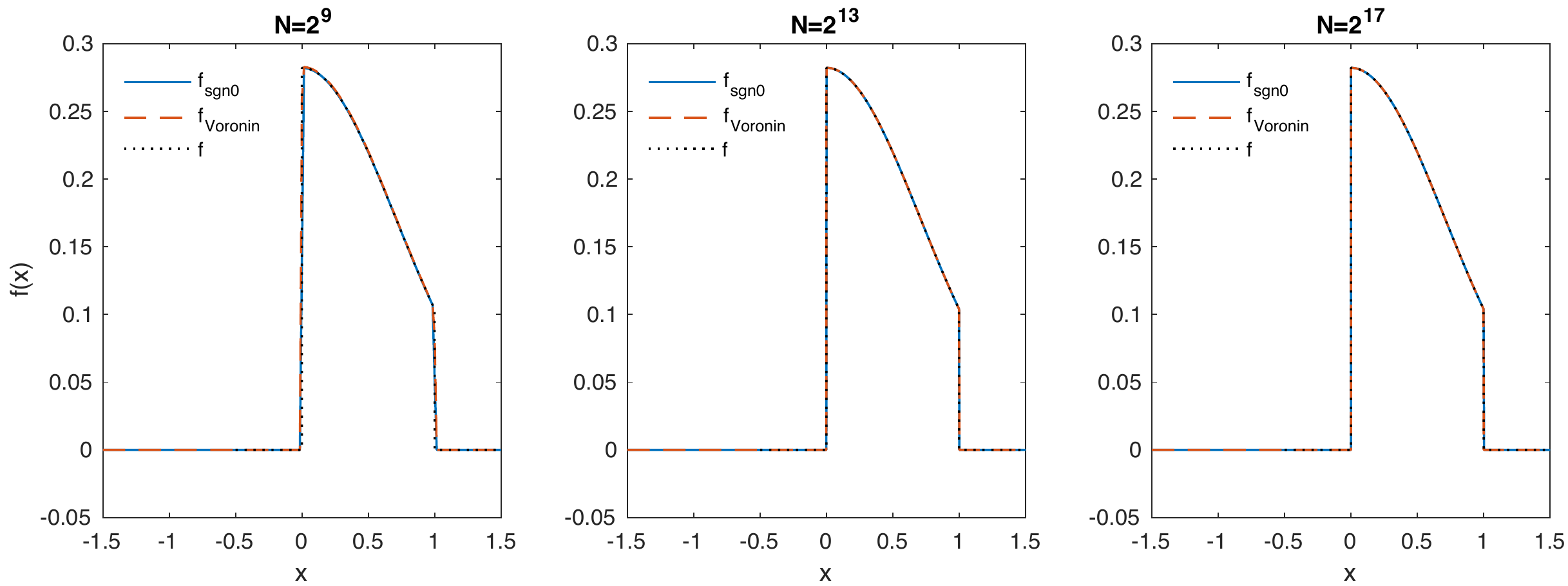}
\caption{Numerical and analytical $f(x)$ using the FFT based method with a symmetrical sign function.}
\label{fig:gausssin0full}
\end{center}
\end{figure}

\begin{figure}[h]
\vspace{-4mm}
\begin{center}
\includegraphics[width=\textwidth]{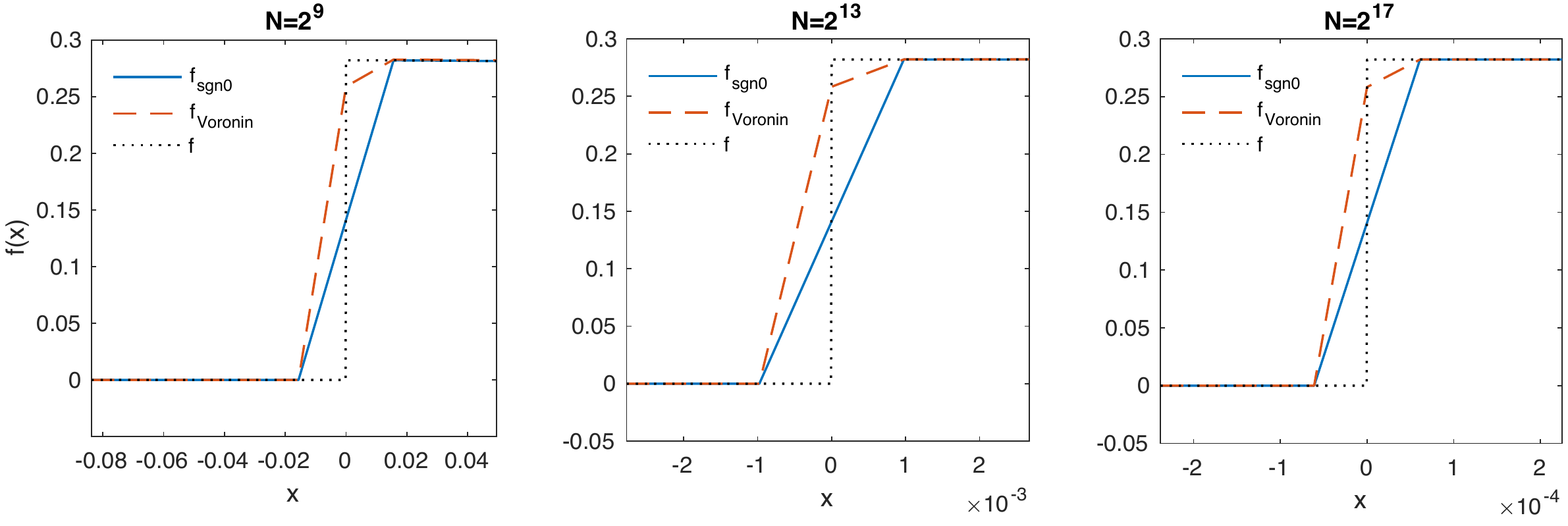}
\caption{Numerical and analytical $f(x)$ using the FFT based method with a symmetrical sign function, focussing on the discontinuity at $x=0$.}
\label{fig:gausssin0up}
\end{center}
\vspace{-4mm}
\end{figure}

\begin{figure}[h]
\vspace{-3mm}
\begin{center}
\includegraphics[width=\textwidth]{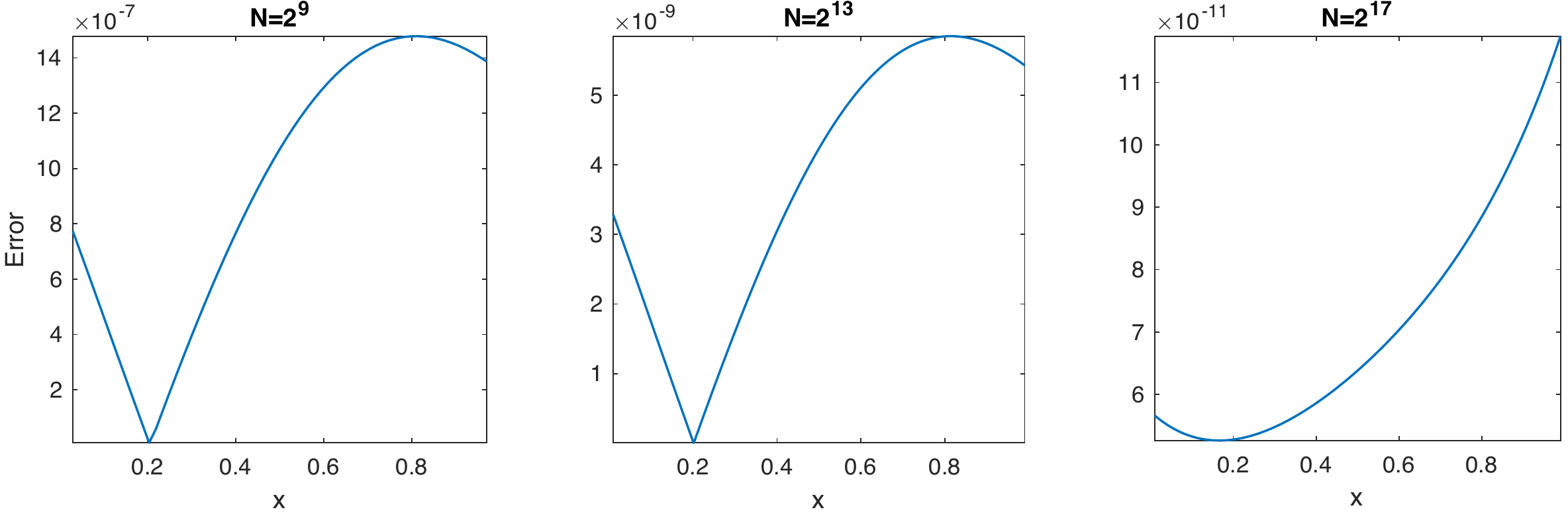}
\caption{Error of the numerical calculation of $f(x)$ for the new iterative Wiener-Hopf method using the sign-based fast Hilbert transform. The scale has been chosen to display the error away from the discontinuities of $f(x)$. The error is calculated by comparing the numerical calculation to the analytic solution.}
\label{fig:gausssin0error}
\end{center}
\vspace{-4mm}
\end{figure}

\begin{figure}[h]
\vspace{-3mm}
\begin{center}
\includegraphics[width=\textwidth]{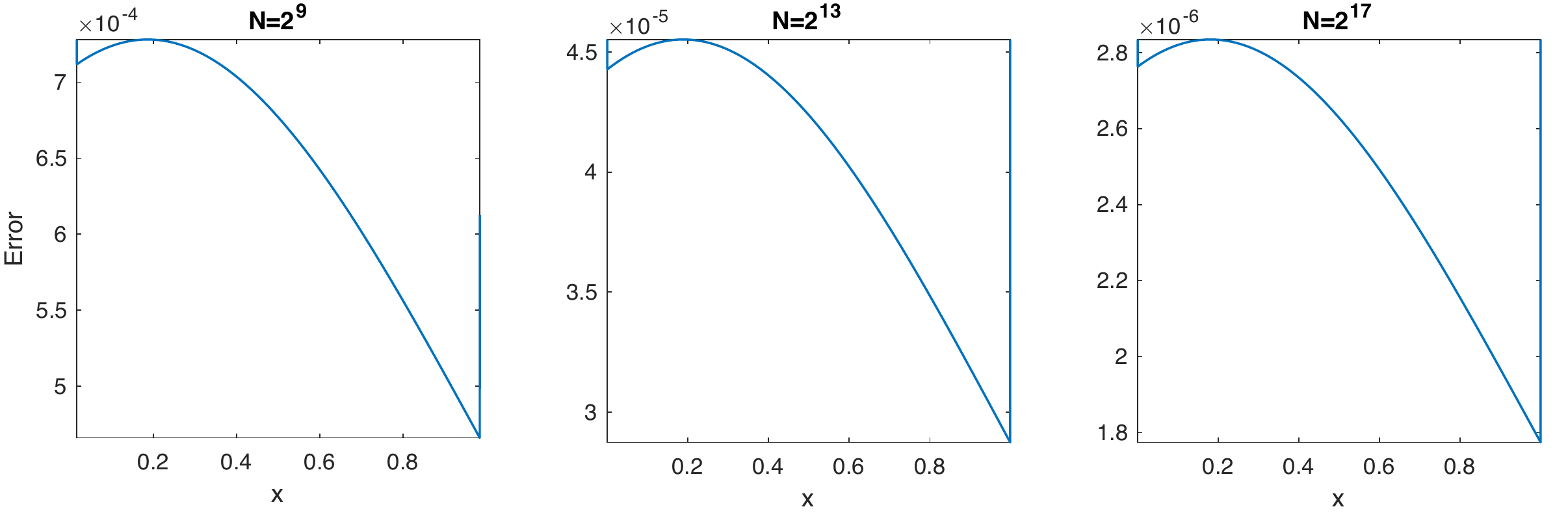}
\caption{Error of the numerical calculation of $f(x)$ for the Voronin method using the sign-based fast Hilbert transform. The scale has been chosen to display the error away from the discontinuities of $f(x)$. The error is calculated by comparing the numerical calculation to the analytic solution.}
\label{fig:gausssin0errorvorr}
\end{center}
\vspace{-4mm}
\end{figure}

The results for the Gaussian test case are shown in Figures \ref{fig:gausssin0full}--\ref{fig:gausssin0error}; $f$ is the analytic solution, $f_{\mathrm{sgn}0}$ and $f_{\mathrm{Vor}}$ are the numerical solutions obtained with the Wiener-Hopf iterative method using the fast Hilbert transform implemented with the symmetric sign function, the latter in the Voronin variant. It is immediately apparent from Figure \ref{fig:gausssin0full} that neither implementation suffers from the overshoot that was seen using the sinc-based fast Hilbert transform. However, looking at the discontinuity more closely in Figure \ref{fig:gausssin0up}, we can see that we will have a peak error at a single state-space grid point as the numerical solution increases to the final value of $f(x)$ more slowly than the analytic function. However, unlike the sinc-based function, where the extent of the oscillations depends not only on the filter, but also the shape the function used, we can state here that as long as the value of $x$ is at least one grid step away from the discontinuity, the answer will be unaffected by the peak error. It is also interesting to note that the error is symmetrical around the discontinuity when the iterative Wiener-Hopf method is used, but not with the Voronin variant. This difference is likely to account for the better performance seen in Figure \ref{fig:gausssin0error} compared to Figure \ref{fig:gausssin0errorvorr}. These display the error results away from the discontinuity and we can see that, although there is some variation in the error across $x$, the results for both methods are superior to those for the sinc-based fast Hilbert transform.

\FloatBarrier

\begin{figure}[h]
\vspace{-2mm}
\begin{center}
\includegraphics[width=\textwidth]{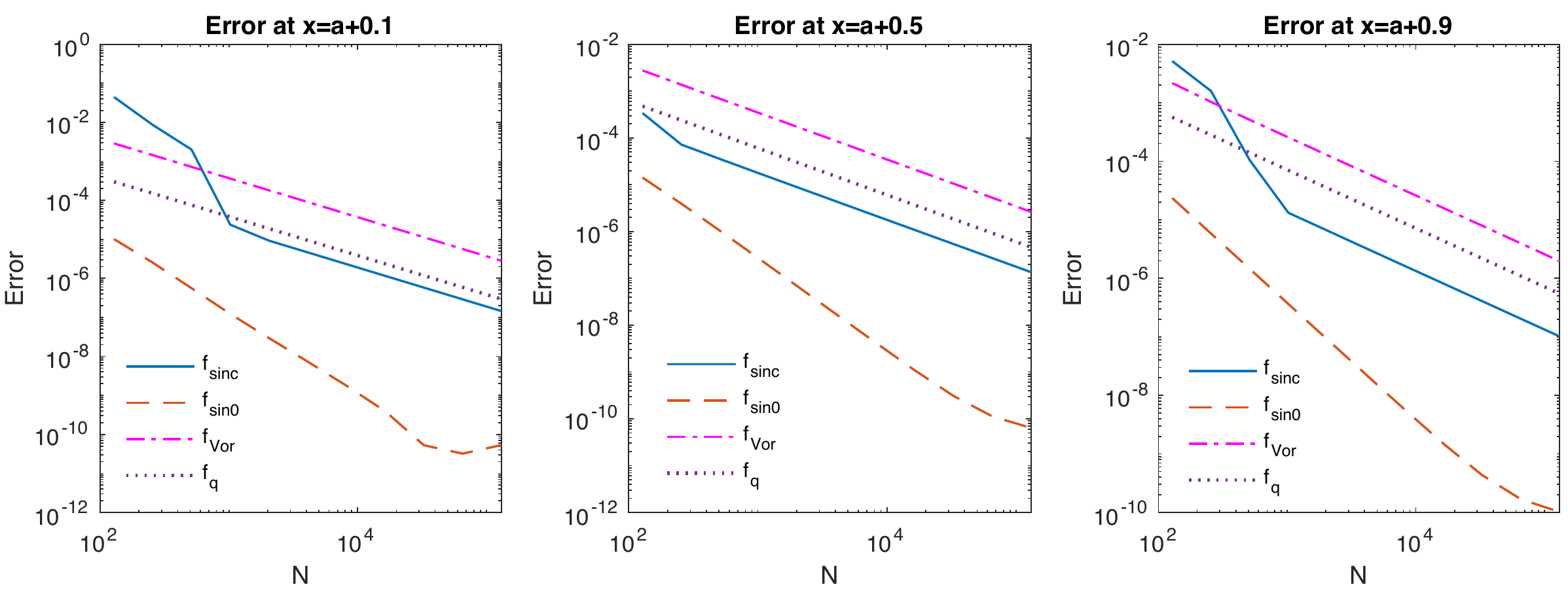}
\caption{Error convergence of the numerical methods vs.\ $N$ with the Gaussian test case.}
\label{fig:gaussconvM}
\end{center}
\vspace{-2mm}
\end{figure}

\begin{figure}[h]
\vspace{-2mm}
\begin{center}
\includegraphics[width=\textwidth]{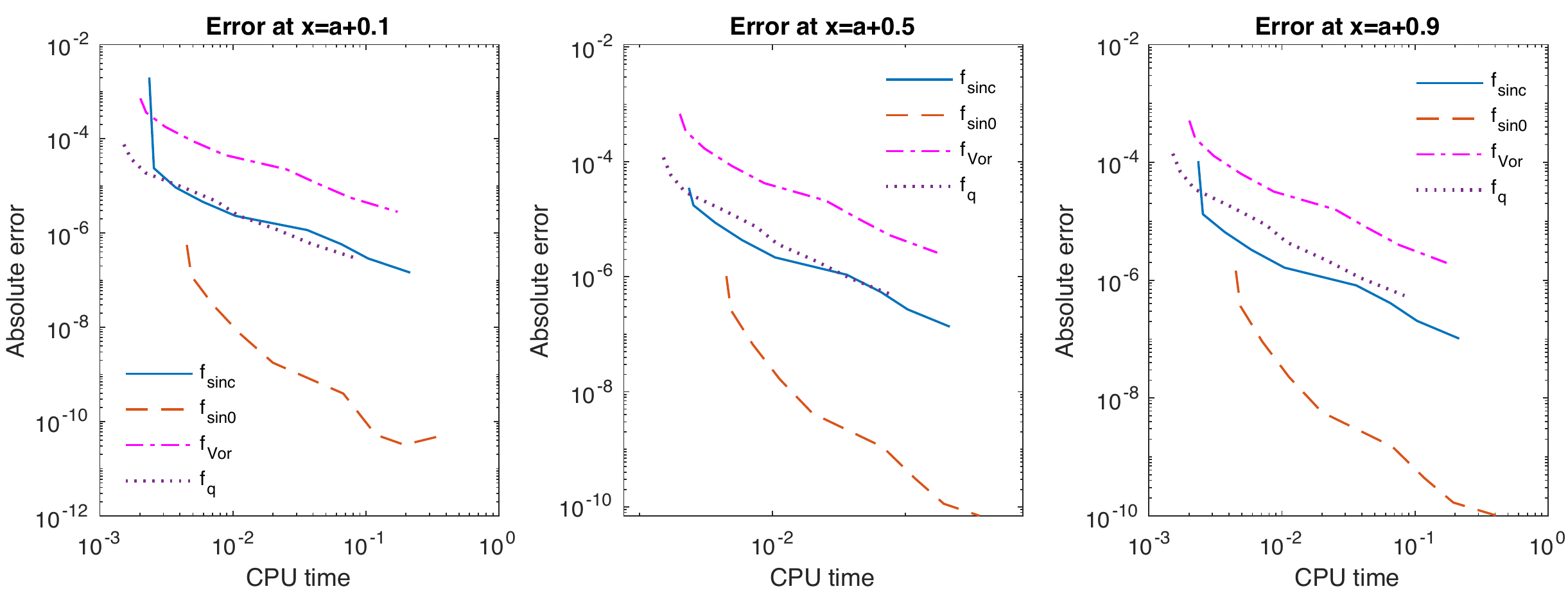}
\caption{Error convergence of the numerical methods vs.\ CPU time with the Gaussian test case.}
\label{fig:gaussconvCPU}
\end{center}
\vspace{-2mm}
\end{figure}

Although it is important to observe the functions which are calculated numerically, when assessing the performance of the numerical methods, the error convergence with CPU time and number of grid points $N$ is also important. We measured this at 10\%, 50\% and 90\% of the range between $a$ and $b$; results for the Gaussian test case are shown in Figures \ref{fig:gaussconvM}--\ref{fig:gaussconvCPU}.% and \ref{fig:cauchconvM}--\ref{fig:expconvCPU}.

The fastest converging method is the Wiener-Hopf iterative method using the sign-based fast Hilbert transform, achieving an error of $O(1/N^2)$. The other methods exhibited $O(1/N)$ error convergence, with the method using the sinc-based fast Hilbert transform with spectral filter achieving better absolute error performance vs.\ $N$ but converging with CPU time almost identically to the quadrature method. The $O(1/N)$ convergence achieved with the sinc-based fast Hilbert transform is consistent with the error bound described by \cite{Stenger1993} for a function with a first order discontinuity, while the $O(1/N^2)$ convergence seen for the sign-based variant is consistent with that reported by \cite{Fusai2016}.

\FloatBarrier
\subsection{Results for Cauchy and Laplace test cases}

\begin{figure}[h]
\vspace{-4mm}
\begin{center}
\includegraphics[width=\textwidth]{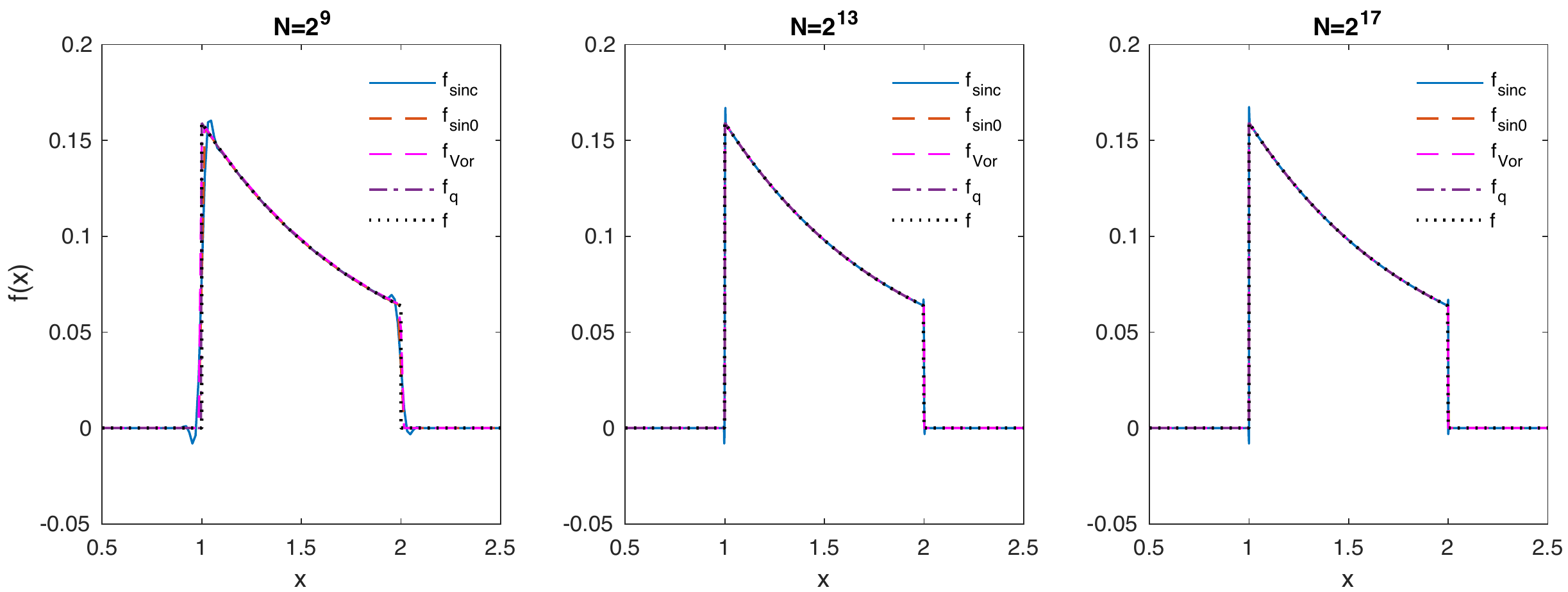}
\caption{Numerical and analytical $f(x)$ with the Cauchy test case.}
\label{fig:cauchallfull}
\end{center}
\vspace{-4mm}
\end{figure}

\begin{figure}[h]
\vspace{-4mm}
\begin{center}
\includegraphics[width=\textwidth]{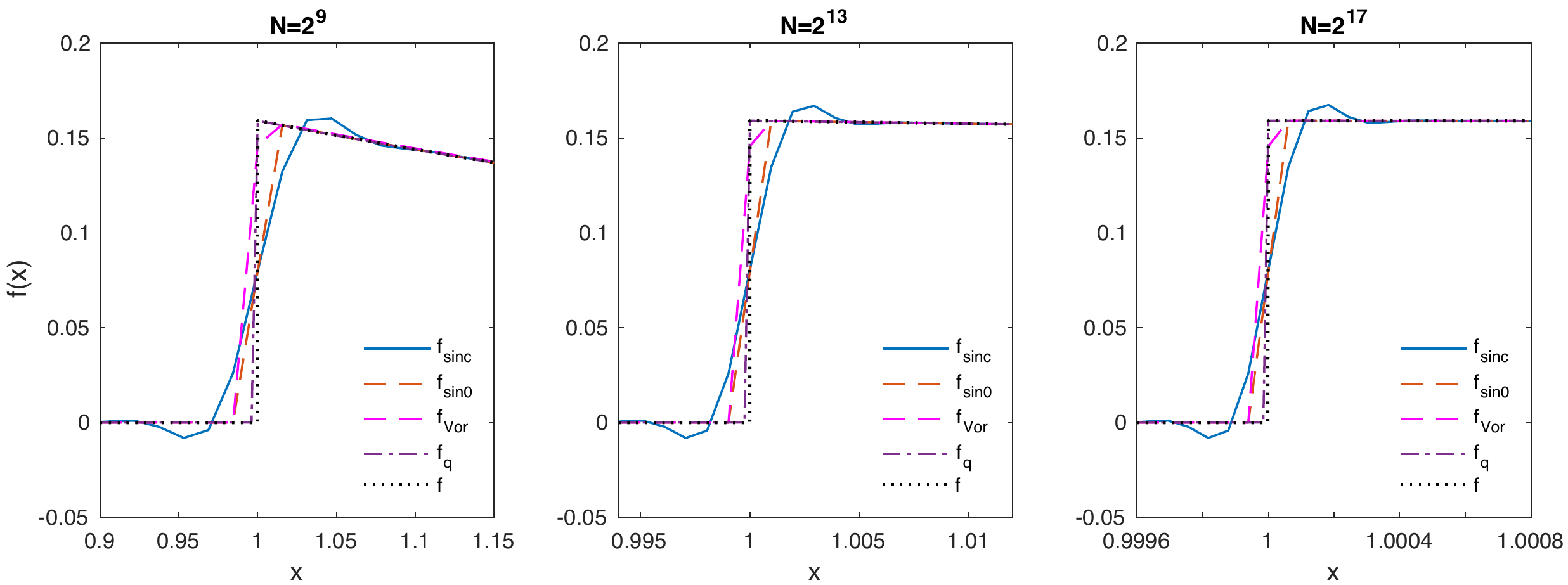}
\caption{Numerical and analytical $f(x)$ with the Cauchy test case focussing on the first discontinuity.}
\label{fig:cauchallup}
\end{center}
\vspace{-4mm}
\end{figure}

\begin{figure}[h]
\vspace{-4mm}
\begin{center}
\includegraphics[width=\textwidth]{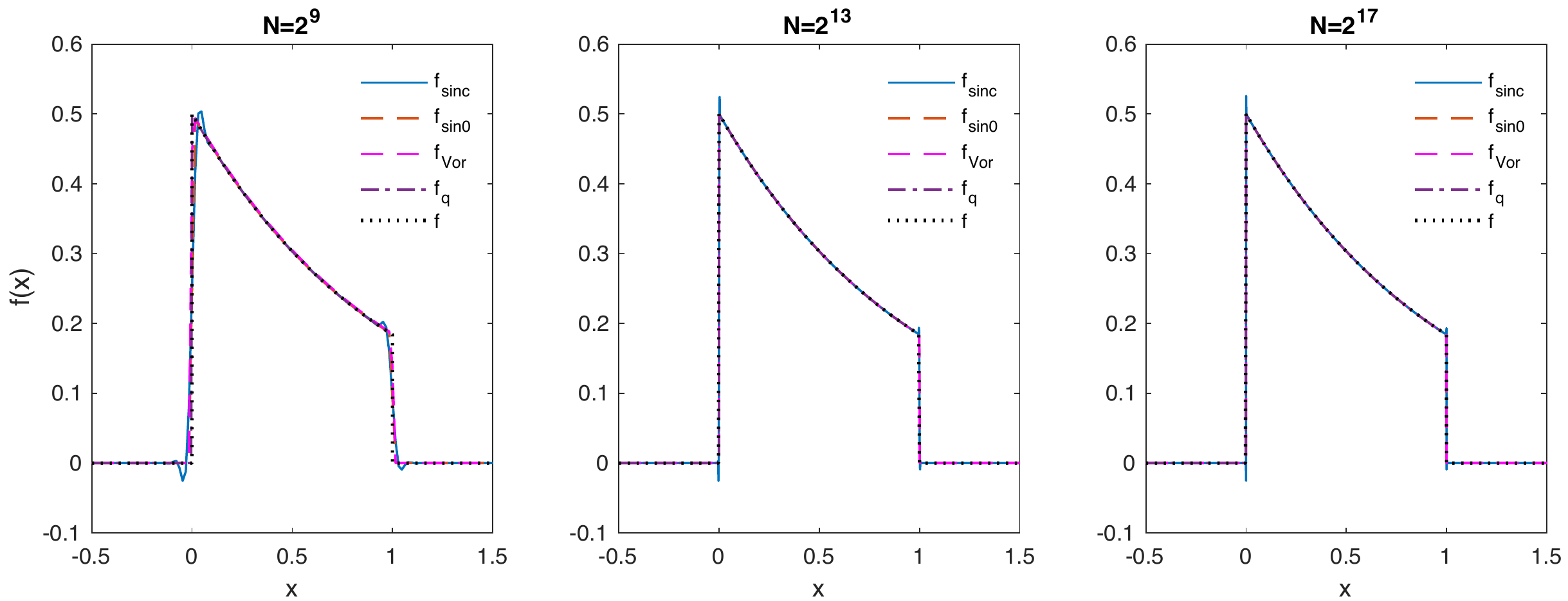}
\caption{Numerical and analytical $f(x)$ with the Laplace test case.}
\label{fig:expallfull}
\end{center}
\vspace{-4mm}
\end{figure}

\begin{figure}[h]
\vspace{-4mm}
\begin{center}
\includegraphics[width=\textwidth]{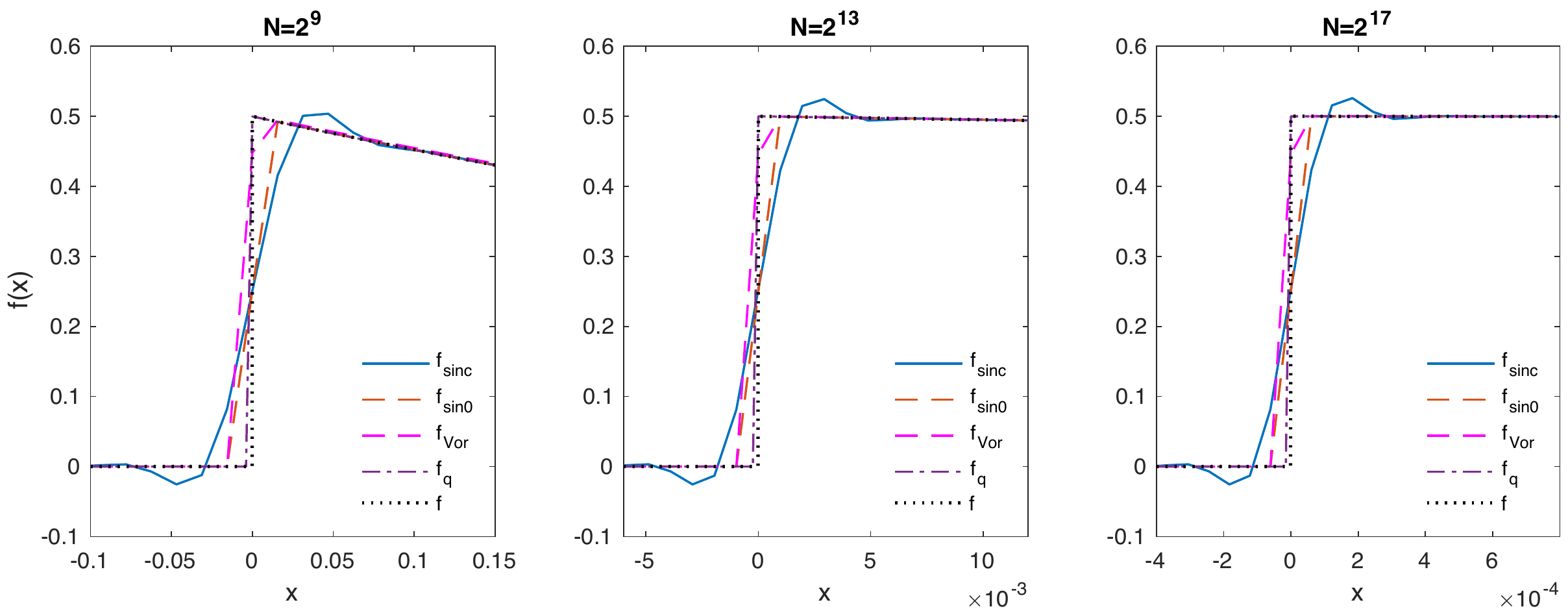}
\caption{Numerical and analytical $f(x)$ with the Laplace test case focussing on the first discontinuity.}
\label{fig:expallup}
\end{center}
\vspace{-4mm}
\end{figure}

Figures \ref{fig:cauchallfull}--\ref{fig:expallup} compare the results for the test cases in Sections \ref{sec:8_Test_Cauch} and \ref{sec:8_Test_Exp} for the iterative Wiener-Hopf method with the sinc- and sign-based fast Hilbert transform methods; in the figures these are labelled $f_\mathrm{sinc}$ and $f_\mathrm{sgn}$. We also compare the performance of the iterative Voronin method with the symmetrical sign function, labelled $f_\mathrm{Vor}$.
An 8th order exponential filter was used with the sinc-based fast Hilbert transform to counteract the oscillations, as described in Section \ref{sec:8_sincfht}. As a benchmark we include results from 4th order quadrature with preconditioner $f_\mathrm{q}$ \citep{Fusai2012,Press2007}, labelled $f_\mathrm{q}$, which was the previous state of the art. The results for the Cauchy and Laplace test cases are consistent with those for the Gaussian test case; the use of the sinc-based fast Hilbert transform results in an overshoot at the function discontinuities and the sign-based method results in a spot error at function discontinuities. We also notice that the quadrature method has a spot error at the discontinuity, but this effects a smaller range of $x$ than our new numerical methods. The reason for this smaller range is that the Fourier-based methods need a truncation in state space at $\pm4(b-a)$ in order to avoid wrap-round effects. In contrast, the range of $x$ for quadrature only needs truncation at the integration limits $a$ and $b$.

\begin{figure}[h]
\vspace{-4mm}
\begin{center}
\includegraphics[width=\textwidth]{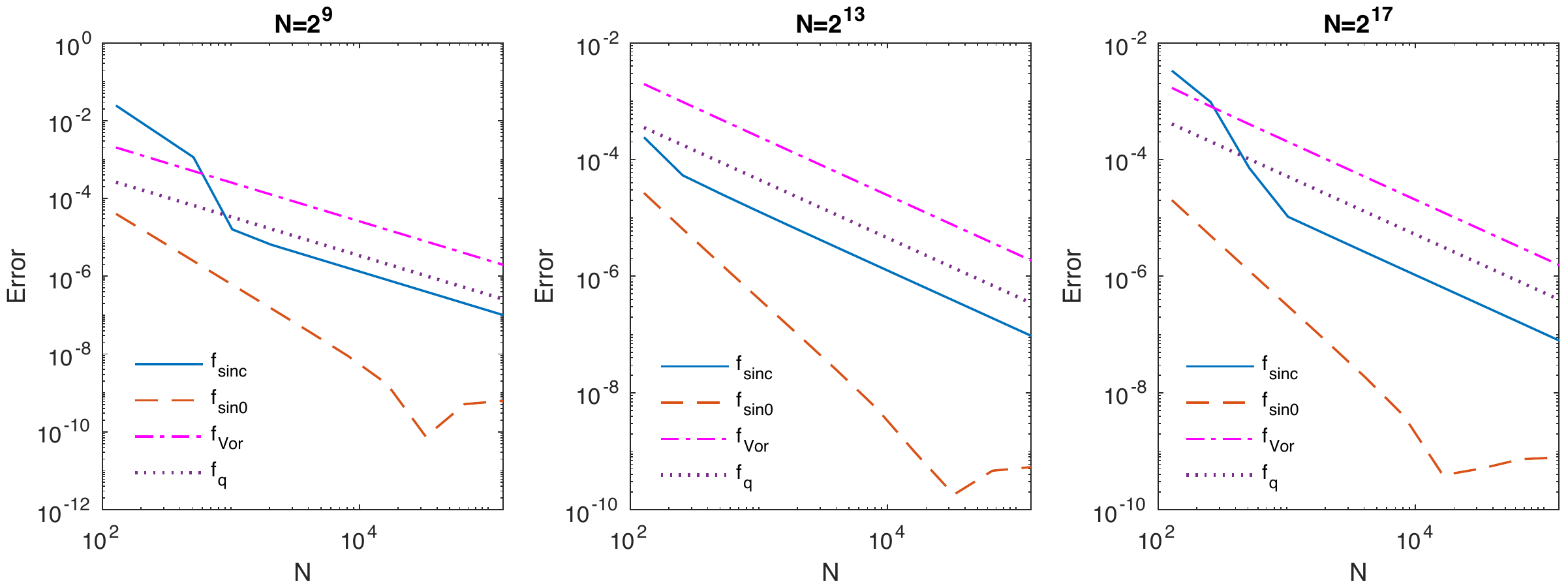}
\caption{Error convergence of the numerical methods vs.\ $N$ with the Cauchy test case.}
\label{fig:cauchconvM}
\end{center}
\vspace{-4mm}
\end{figure}

\begin{figure}[h]
\vspace{-4mm}
\begin{center}
\includegraphics[width=\textwidth]{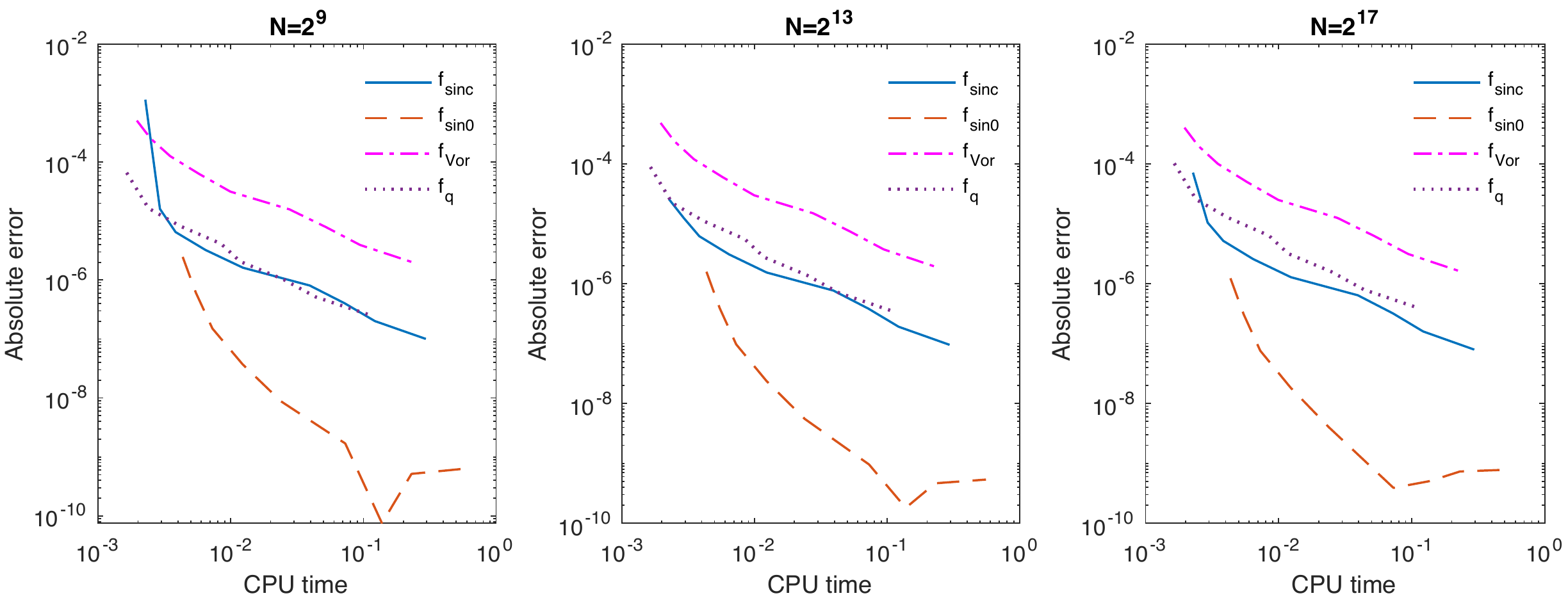}
\caption{Error convergence of the numerical methods vs.\ CPU time with the Cauchy test case.}
\label{fig:cauchconvCPU}
\end{center}
\vspace{-4mm}
\end{figure}

\begin{figure}[h]
\vspace{-4mm}
\begin{center}
\includegraphics[width=\textwidth]{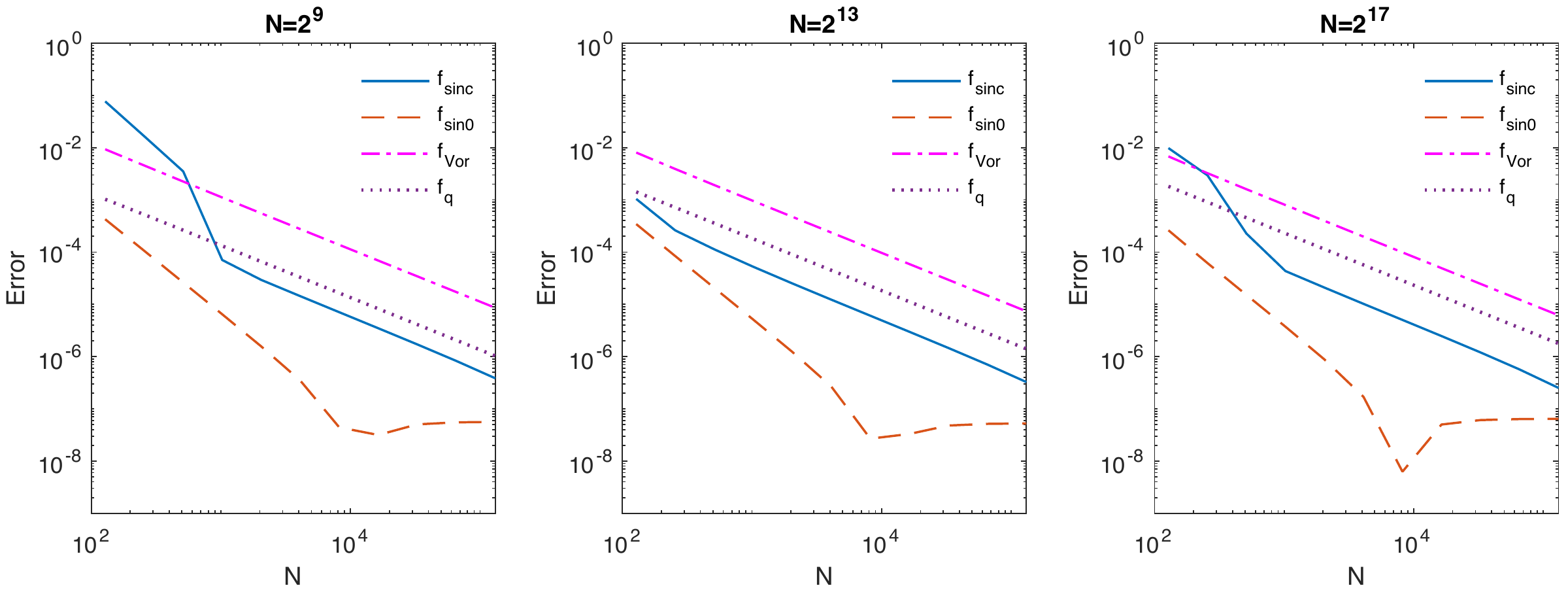}
\caption{Error convergence of the numerical methods vs.\ $N$ with the Laplace test case.}
\label{fig:expconvM}
\end{center}
\vspace{-4mm}
\end{figure}

\begin{figure}[h]
\vspace{-4mm}
\begin{center}
\includegraphics[width=\textwidth]{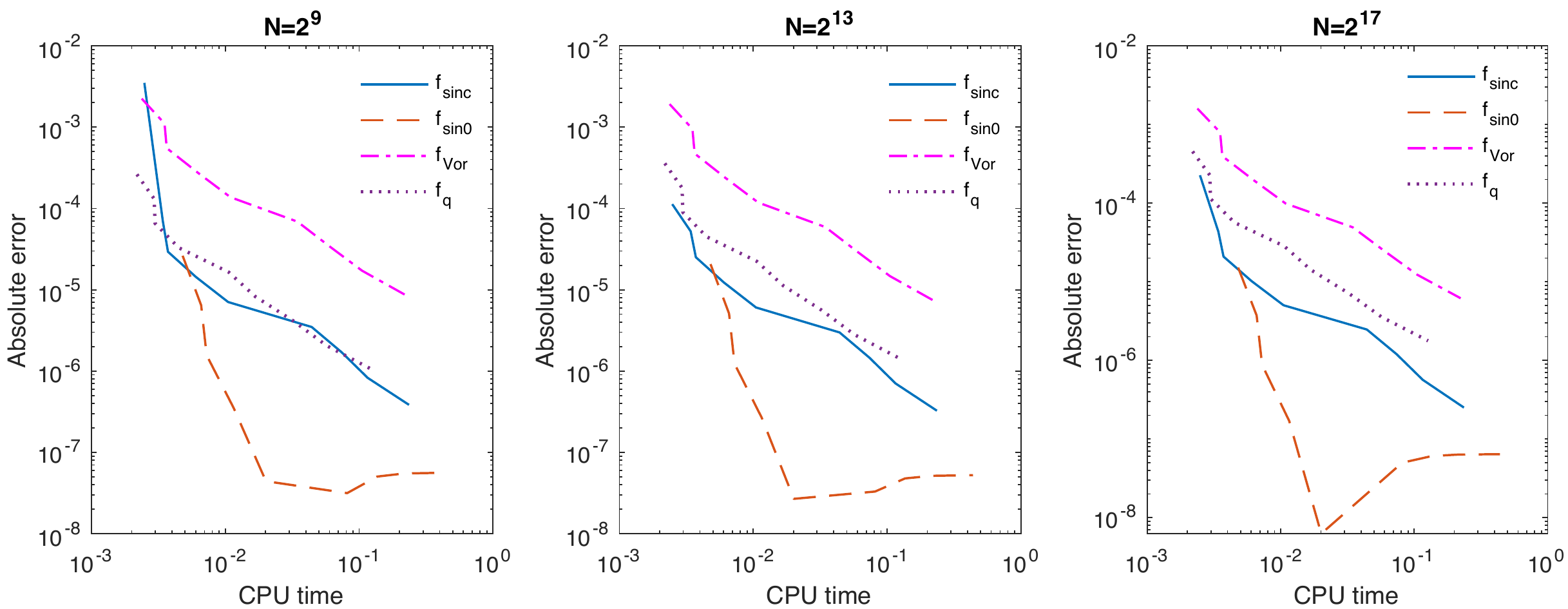}
\caption{Error convergence of the numerical methods vs.\ CPU time with the Laplace test case.}
\label{fig:expconvCPU}
\end{center}
\vspace{-4mm}
\end{figure}

We also measured the error convergence with the Cauchy and Laplace test cases and the results are shown in Figures \ref{fig:cauchconvM}--\ref{fig:expconvCPU}. These confirm the findings with the Gaussian test case in Section \ref{sec:8_Test_Gauss} which showed that the best performing method is the new iterative solution to the Wiener-Hopf equation with the Hilbert transform implemented using the FFT with the symmetrical sign function.

\FloatBarrier

\section{Conclusion}
We implemented four methods to solve general Fredholm equations of the second kind and assessed their performance for three test cases with analytical solutions. The methods are: 4th order quadrature with preconditioner \cite{Fusai2012,Press2007}; two iterative solutions that use the Wiener-Hopf method with the sinc- or sign-based fast Hilbert transform; a variant of the iterative method based on Voronin's partial solution to the Fredholm equation with the sign-based fast Hilbert transform.

Unlike an earlier application in option pricing with exponential convergence of a weighted average error, the iterative Wiener-Hopf method with the sinc-based fast Hilbert transform does not turn out optimal for a general solution of the Fredholm equation, having $O(1/N)$ convergence with the number of FFT grid points $N$ and high errors close to the function discontinuities; this can be explained with the different requirements of the two problems, as here the solution over the whole interval is required. Instead, the iterative Wiener-Hopf method with the sign-based fast Hilbert transform has $O(1/N^2)$ convergence and therefore performs better than its sinc-based sibling, the quadrature method from the literature and the iterative method based on Voronin's partial solution, whose convergence is $O(1/N)$ even with the sign-based fast Hilbert transform. So in terms of error convergence, the iterative Wiener-Hopf method with the sign-based fast Hilbert transform reveals the new state of the art for the numerical solution of general Fredholm equations, achieving double the convergence speed of the known 4th order quadrature method.

The other aspect which we must compare for the different methods is the peak error at a discontinuity of $f(x)$ as shown in Figures \ref{fig:expallup} and \ref{fig:cauchallup}. This error is wider for the Wiener-Hopf method than for the quadrature method because a wider $x$ range is required to avoid the wrap-around effects of Fourier transform methods. Therefore, if we require an accurate answer close to a discontinuity, the quadrature method may be best. However, the excellent CPU time vs.\ error performance shown in Figures \ref{fig:gaussconvCPU}, \ref{fig:cauchconvCPU} and \ref{fig:expconvCPU} recommend using the Wiener-Hopf method with the sign-based fast Hilbert transform and a larger grid size to yield the required accuracy close to the discontinuities.

%\bigskip
%\textcolor{red}{\textit{
%To do:
%\begin{enumerate}
%\item Perhaps discuss slightly more in depth the iterative solution by Dionisios Margetis and Jaehyuk Choi, Studies in Applied Mathematics \textbf{117}, 1--25 (2006) mentioned at the start of Section 4.3.
%\item Look at the s of Dean G.~Duffy for examples of WH/F equations not from finance, more Fredholm equations with analytical solution, and WHE with continuous forcing function and analytical solution.
%\item In the iterative method, monitor not only the behaviour and the error of $\widehat{f_0}$, but also of the auxiliary variables which actually enter the iteration.
%\item Obtain a proof of $O(1/n^2)$ performance.
%\item Optimise filtering to balance overshoot and error, look at Planck taper as well as exponential filter. Briefly mention that this does not bring an improvement briefly mention in a sentence.
%\item Check the analysis of the differences between this method and the option pricing method in order to understand why the sinc-based method only achieves $O(1/N)$ convergence: the 1/n performance is explained referencing the error bound described by Stenger for a first order discontinuity.
%\item Analyse why addition of zero padding degrades performance.
%\item Perhaps add some detail to the description of the quadrature method in the first paragraph of Section 2.4.
%\end{enumerate}}}

\section*{Acknowledgments}

We are thankful to Carlo Sgarra for a cue about Sec.~\ref{sec:Noble}, and to Vito Daniele for useful discussions on the difficulties of factorising the matrices $\mathbf{L}$ of Sec.~\ref{sec:Noble} and $\mathbf{M}$ of Sec.~\ref{sec:Voronin}, as well as for a preprint of his book a year before it was published \citep{Daniele2014}. % and for reading a preprint of this paper.

The support of the Economic and Social Research Council (ESRC) in funding the Systemic Risk Centre (grant number ES/K002309/1) and of the Engineering and Physical Sciences Research Council (EPSRC) in funding the UK Centre for Doctoral Training in Financial Computing and Analytics (grant number 1482817) are gratefully acknowledged.

\bibliography{paper}
%Explore the biblatex package
%\bibliographystyle{unsrt}
\bibliographystyle{imamatgeneric}

\end{document}